\documentclass[12pt]{article}
\usepackage{pict2e}
\usepackage{bbm}
\usepackage{amsfonts}
\usepackage{pifont}
\usepackage{hyperref}
\usepackage{mathrsfs}
\usepackage{indentfirst}
\usepackage{amsmath}
\usepackage{amssymb}
\usepackage[amsmath, thmmarks]{ntheorem}
\usepackage[noadjust]{cite}
\usepackage{graphicx}
\usepackage{tikz}
\usepackage{newunicodechar}
\newunicodechar{ }{\leavevmode\,} % U+202F 映射成LaTeX细空格\,
\newtheorem{definition}{\bf Definition}[section]
\newtheorem{lemma}{\bf Lemma}[section]
\newtheorem{theorem}{\bf Theorem}[section]
\newtheorem{remark}{\bf Remark}[section]
\newtheorem{corollary}{\bf Corollary}[section]
\newtheorem{proposition}{\bf Proposition}[section]
\newtheorem{example}{\bf Example}[section]

\def\QEDopen{{\setlength{\fboxsep}{0pt}\setlength{\fboxrule}{0.2pt}\fbox{\rule[0pt]{0pt}{1.3ex}\rule[0pt]{1.3ex}{0pt}}}}
\def\QED{\QEDopen}
\def\proof{{\bf Proof.} }
\def\endproof{\hspace*{\fill}~\QED\par\endtrivlist\unskip}
\begin{document}
\setcounter{page}{1}

\title{{\textbf{Finite groups with consistent subgroup lattices}}
\thanks {Supported by the National Natural Science Foundation of China (Nos. 12471440 and 12671562)}
}
\author{Peng He\footnote{\emph{E-mail address}: 443966297@qq.com}, Xue-ping Wang\footnote{Corresponding author. xpwang1@hotmail.com; fax: +86-28-84761502}\\
\emph{\small (School of Mathematical Sciences, Sichuan Normal University}\\ \emph{\small Chengdu 610066, Sichuan, People's Republic of China)}}
\newcommand{\pp}[2]{\frac{\partial #1}{\partial #2}}
\date{}
\maketitle
\begin{quote}
{\bf Abstract:} {Motivated by the problem of Gragg and Kung on characterizing finite groups with consistent subgroup lattices, we investigate 
how the consistency of the subgroup lattice $L(G)$ constrains the structure of the underlying finite group $G$. To this end, we first establish 
several equivalent characterizations of consistency for lattices of finite length in terms of operation-preserving subposets and forbidden sublattices. 
More interestingly, for complete atomistic lattices of finite length, we introduce the consistency exchange property (CEP) for the canonical closure operator and show that 
consistency is equivalent to CEP, and that CEP is strictly weaker than the Steinitz–Mac Lane exchange property. We then apply this framework 
to semidirect products $G=V\rtimes_{\alpha} C_m$, where $V$ is an elementary abelian $p$-group, $m$ is square-free and $p\nmid m$. Writing $d=o(\alpha)$ and
letting $C_\alpha=\langle\alpha\rangle$, we prove that $L(G)$ is consistent if
and only if either $d=1$, or $d$ is prime and $V$ is isotypic as an
$\mathbb F_pC_\alpha$-module. For prime $d$, this latter condition is
equivalent both to the irreducibility of the minimal polynomial of $\alpha$
over $\mathbb F_p$ and to $\mathbb F_p[\alpha]$ being a field.
More generally, we verify that if $G=V\rtimes_\theta H$ where $H$ acts faithfully on $V$ and $p\nmid |H|$, then $L(G)$ is consistent if and only if either $H$ is trivial, or $H$ is cyclic of prime order and $V$ is isotypic as an $\mathbb F_pH$-module. We thus obtain a complete solution to the problem of Gragg and Kung for this class of coprime semidirect products.}

\end{quote}
\section{Introduction}
For a finite group $G$, the subgroup lattice $L(G)$ consists 
of all subgroups of $G$, ordered by inclusion, with meet, $\wedge$, 
given by intersection and join, $\vee$, given by subgroup
generation \cite{Suz, Schmidt}. 
A fundamental problem is to determine how much information about $G$ 
is encoded in $L(G)$. The investigation of this problem began 
with the work of Rottl\"{a}nder \cite{Rott}. Since then, 
many authors established a series of characterization theorems 
for finite groups using classical lattice properties including 
modularity and distributivity. One of the best-known results is 
Ore's characterization of locally cyclic groups in terms of 
the distributivity of their subgroup lattices \cite{Ore}.  
The systematic study of the relationship between groups and 
their subgroup lattices was subsequently developed by 
Suzuki \cite{Suz51,Suz52}. Subsequent work has further 
extended this framework. For example, Andreeva, Schmidt, 
and Toborg \cite{Andreeva} investigated lattice-defined classes of finite 
groups with modular Sylow subgroups, Shareshian and 
Woodroofe \cite{Shareshian} obtained a subgroup lattice characterization 
of finite solvable groups, and Jovanovi\'{c}, \v{S}e\v{s}elja, 
and Tepav\v{c}evi\'{c} \cite{Jovanovi,Jovanovi24} studied lattice 
characterizations of finite nilpotent groups. 

More recently, modularity 
in the subgroup lattices of finite soluble groups and of polycyclic groups was investigated in \cite{Liu} and \cite{Trombetti}, respectively. 
Related questions concerning prescribed intervals in subgroup 
lattices were examined by Baddeley and Lucchini \cite{Baddeley} 
and by Aschbacher \cite{Aschbacher}, while combinatorial structures 
arising from the notion of independence 
in finite groups were studied in \cite{Lucchini}. These results demonstrate that classical modularity and distributivity 
impose strong restrictions on finite groups and provide powerful tools 
for characterizing their structure. The consistency of a lattice was 
introduced by Kung \cite{Kung} when he studied matchings between 
join-irreducible and meet-irreducible elements of finite lattices. 
Since consistency is a natural extension of modularity, Gragg and Kung 
posed the following problem in \cite[Problem 8.2]{Gragg}: 
\textit{Characterize the finite groups with a consistent subgroup lattice.} 
More precisely, they asked which finite groups $G$ have the property 
that $L(G)$ is a consistent lattice. 

To address the above problem,  we first establish structural characterizations of consistency for lattices of finite length. 

Let $L$ be a lattice. Recall that an element $j\in L$ is join-irreducible 
if $j=a\vee b$ implies that $j=a$ or $j=b$. We write  $J(L)$ for the 
set of all join-irreducible elements of $L$ and $J^{+}(L)=J(L)\setminus\{0\}$, where $0$ is the least element of $L$. 
A join-irreducible element $j\in J(L)$ is called consistent if $x\vee j\in J([x,1])$ for every $x\in L$. 
A lattice of finite length is 
called consistent if every element of $J(L)$ is consistent. Since $0$ 
is automatically consistent, it is enough to consider the elements 
of $J^{+}(L)$ \cite{Stern99}. Following Kung's introduction of consistency, 
further structural properties of consistent lattices were investigated by Walendziak \cite{Walendziak94, Walendziak95} and others; 
see also \cite{Stern99}. In particular, Gragg and Kung \cite{Gragg} 
studied consistent dually semimodular lattices and established several 
connections between consistency, atomicity and modularity. 
Moreover, Stern characterized consistency in semimodular 
lattices of finite length by means of special forbidden semimodular 
sublattice $S_7$ \cite[Theorem 6.5.1]{Stern99}. 
Accordingly, it is natural to ask whether consistent lattices of finite 
length can be characterized by forbidden substructures. 
The first main result of this article, presented below 
as Theorem \ref{TT01}, provides a complete answer to this question.
\begin{theorem}\label{TT01}
Let $L$ be a lattice of finite length. Then the following conditions are equivalent:\\
(1) $L$ is consistent;\\
(2) $L$ contains no operation-preserving subposet of the form shown in Figure 1, with the distinguished element $j\in J^{+}(L)$; \\
(3) $L$ contains no sublattice of any of the three 
forms $L_5, S_7$ and $B^{*}_{3}$ shown in Figure 4, 
with the distinguished element $j\in J^{+}(L)$.   
\end{theorem}

On the other hand, Reuter \cite{Reuter} proved that a finite lattice 
is consistent if and only if it satisfies the Kurosh-Ore 
replacement property. Walendziak \cite{Walendziak95} subsequently 
provided a closure-theoretic characterization of consistency: 
a finite lattice is the lattice of closed sets of a closure 
space with the Steinitz exchange property if and only if 
it is consistent. It is well known that every geometric lattice 
is a finite atomistic lattice and admits an equivalent representation 
by a closure operator satisfying the Steinitz-Mac Lane exchange property. 
Geometric lattices are consistent, but consistent lattices are not 
necessarily geometric \cite{Stern99}. This raises the question of 
whether the consistency of atomistic lattices admits an equivalent 
characterization in terms of a closure operator satisfying 
another exchange property. To formulate and answer this question precisely, 
we first recall the relevant definitions and properties as follows. 

Let $L$ be a lattice. For $x, y\in L$, we say that $y$ covers $x$, 
and write $y\succ x$, if $x< y$ and the interval 
$$[x, y]=\{z\in L: x\leq z\leq y\}=\{x,y\}.$$ 
An element that covers the least element $0$ of $L$ is 
called an atom of $L$. The set of all atoms in $L$ is denoted by $A(L)$. 
Moreover, for each $x\in L$, we define $A(x)=\{a\in A(L): a\leq x\}$.
A complete lattice $L$ is called atomistic 
if $$x=\bigvee\{a\in A(L): a\leq x\}=\bigvee A(x)$$ 
for every $x\in L$ \cite{Crawley73}. For a complete atomistic 
lattice $L$ of finite length, 
one has $A(L)=J^{+}(L)$ \cite{Birkhoff73, Davey}. 
Let $L$ be a complete lattice of finite length. 
Define the map $\mbox{cl}: 2^{J^{+}(L)}\longrightarrow 2^{J^{+}(L)}$ by
\begin{equation}\label{eqqq-1}
\mbox{cl}(X)=\{x\in J^{+}(L): x\leq \bigvee X\}
\end{equation}
for every $X\subseteq J^{+}(L)$. Then $\mbox{cl}$ is a closure 
operator on $J^{+}(L)$ and $L\cong \mathcal{C}(J^{+}(L), \mbox{cl})$, 
where $\mathcal{C}(J^{+}(L), \mbox{cl})$ denotes the complete lattice 
of all closed subsets of $J^{+}(L)$ ordered by 
inclusion \cite{Nation, Davey}. The closure operator $\mbox{cl}$ defined by \eqref{eqqq-1} is called the canonical closure operator associated with $L$. 

Recall that a closure operator 
$\mbox{cl}: 2^I\longrightarrow 2^I$ satisfies the Steinitz–Mac Lane exchange property (EP) if, for any $x,y \in I$ and $X\subseteq I$, 
$$x\in \mbox{cl}(X\cup\{y\}) \setminus \mbox{cl}(X)\Longrightarrow y\in \mbox{cl}(X\cup\{x\}) \setminus \mbox{cl}(X).$$ Suppose now that $L$ is a 
complete atomistic lattice and that $\mbox{cl}$ is the 
canonical closure operator associated with $L$.
Then $L\cong \mathcal{C}(A(L), \mbox{cl})$ since $J^{+}(L)=A(L)$. 
If, in addition, $\mbox{cl}$ satisfies EP, 
then the lattice $\mathcal{C}(A(L), \mbox{cl})$ is a complete 
atomistic semimodular lattice \cite{Crawley73}. Moreover, if $L$ is finite, 
then $L$ is a geometric lattice \cite{Birkhoff73}. 
Motivated by the classical Steinitz-Mac Lane exchange property, 
we introduce the consistency exchange property as follows.
\begin{definition}\label{WEK}
\emph{Let $\mbox{cl}$ be a closure operator on a nonempty set $I$. We say that $\mbox{cl}$ satisfies the consistency exchange property (CEP) if, for any $X\subseteq I$ and $x,y,z\in I$,
$$y\in \mbox{cl}(X\cup\{x,z\}) \setminus [\mbox{cl}(X\cup\{x\})\cup \mbox{cl}(X\cup \{z\})]$$
implies $$\mbox{cl}(X\cup\{x,z\})\neq \mbox{cl}(X\cup\{y\}).$$}
\end{definition}
We show that EP implies CEP, but the converse fails. 
Recall that a sublattice $L_0$ of a lattice $L$ with the least 
element $0$ is called a $0$-sublattice if $0\in L_0$.  We then obtain equivalent characterizations of consistency 
in complete atomistic lattices of finite length in terms 
of CEP and forbidden $0$-sublattices, respectively. 
\begin{theorem}\label{T04}
Let $L$ be a complete atomistic lattice of finite length and let $\emph{cl}$ denote the canonical closure operator associated with $L$. Then the following conditions are equivalent:\\
(i) $L$ is consistent.\\
(ii) $\emph{cl}$ satisfies CEP.\\
(iii) $L$ contains no $0$-sublattice of the form $L_{5}^{u}$ shown in Figure 8, with the distinguished element $j\in J^{+}(L)$.\\
(iv) $L$ contains no $0$-sublattice of 
the form $L_{5}^{l}$ shown in Figure 9, with the 
distinguished element $j\in J^{+}(L)$.
\end{theorem}

Equipped with the above lattice-theoretic characterizations 
of consistency, we now turn to Gragg and Kung's problem of 
characterizing finite groups with a consistent subgroup lattice.

The subgroup lattice of the alternating group $A_4$ provides an 
instructive example \cite{Stern99}. Indeed, $L(A_4)$ is consistent 
but not modular, and it is neither upper nor 
lower semimodular (see also \cite{Gragg}). 
Moreover, $A_4\cong C^2_{2}\rtimes_{\theta} C_3,$ 
where $C_3$ acts faithfully on the elementary abelian 
group $C^2_{2}$ \cite{Scott, Rotman}. This demonstrates that 
a nonmodular consistent subgroup lattice can be obtained 
from a nontrivial linear action on an elementary abelian group.
Motivated by this observation, 
we investigate semidirect products $$G=V\rtimes_{\theta} H,$$ 
where $V\cong C^n_p$ is an elementary abelian $p$-group, $n\geq 1$, $p$ is prime, and $H$ is a finite group with $p\nmid |H|$. 
Note that, throughout this article, we identify $V\cong C^{n}_{p}$ 
with its natural structure as an 
$n$-dimensional vector space over $\mathbb{F}_p$, write its operation 
additively, and denote its zero element by $0$. 
The action of $H$ on $V$ 
is specified by a homomorphism 
$$\theta: H\longrightarrow \mbox{GL}_{\mathbb{F}_{p}}(V),$$
which is said to be faithful if $\theta$ is injective, equivalently, 
if no nonidentity element of $H$ acts trivially on $V$.
The action of $H$ endows $V$ with the structure of an $\mathbb{F}_{p}H$-module. 
Since $p\nmid |H|$, Maschke's theorem implies that $V$ is semisimple.
Recall that an $\mathbb{F}_{p}H$-module $V$ is called isotypic 
if $V\cong S^{\oplus l}$ for some simple $\mathbb{F}_{p}H$-module $S$ 
and some integer $l\geq 1$ \cite{Lux, Webb}. Equivalently, $V$ is isotypic 
if and only if all its irreducible constituents are mutually 
isomorphic.  Therefore, the main group-theoretic question addressed in this 
article is the following: to what extent does the consistency 
of $L(G)$ determine the structure of the acting group $H$ 
and the corresponding $\mathbb{F}_pH$-module $V$? 

We first investigate cyclic actions. Recall that for $T\in \mbox{End}_{\mathbb{F}_{p}}(V)$, its minimal polynomial $m_T(x)$ is the unique monic polynomial of least positive degree satisfying $m_{T}(T)=0$. 
The $\mathbb{F}_{p}$-subalgebra generated by $T$ is $$\mathbb{F}_{p}[T]=\{f(T):f(x)\in \mathbb{F}_{p}[x]\}\subseteq \mbox{End}_{\mathbb{F}_{p}}(V).$$
Moreover, $$\mathbb{F}_{p}[T]\cong \mathbb{F}_{p}[x]/(m_{T}(x)),$$ so $\mathbb{F}_{p}[T]$ is a field if and only if $m_T(x)$ is irreducible over $\mathbb{F}_{p}$ \cite{Axler}. 
Using Maschke's theorem, the structure of cyclic group algebras and Theorem \ref{T04}, we obtain the following main result.
\begin{theorem}\label{T001}
Let $m$ be square-free, $p \nmid m$, $C_m =\langle b \rangle$, and let $\alpha\in \emph{GL}_{\mathbb F_p}(V)$ satisfy $\alpha^m=\emph{id}_V$. Let $C_m$ act on $V$ 
by $b\cdot v=\alpha(v)$.
Set $G=V\rtimes_{\alpha} C_m$, $o(\alpha)=d$ and $C_\alpha=\langle \alpha \rangle\leq \emph{GL}_{\mathbb F_p}(V)$. Then the following conditions are equivalent.\\
$(a)$ $L(G)$ is consistent.\\
$(b)$ Either $d=1$, or $d$ is prime and $V$ is isotypic as an $\mathbb F_pC_{\alpha}$-module.\\ 
$(c)$ Either $d=1$, or $d$ is prime and $m_{\alpha}(x)$  is irreducible over $\mathbb{F}_p$.\\
$(d)$ Either $d=1$, or $d$ is prime 
and $\mathbb{F}_p[\alpha]$ is a field.
\end{theorem}
We then extend this characterization to arbitrary 
faithful actions of finite groups.

\begin{theorem}\label{TH4.3}
Let $G=V\rtimes_{\theta} H$, where $H$ is a finite group acting faithfully on $V$ and $p\nmid |H|$. Then $L(G)$ is consistent if and only if either $H=\{1\}$, or $H\cong C_q$ for some prime $q$ and $V$ is isotypic as an $\mathbb{F}_pH$-module.
\end{theorem}
These results demonstrate that the consistency of a subgroup lattice imposes strong structural restrictions on the underlying group and its linear actions.

We use the following notations throughout this article. Unless stated otherwise, all groups considered below are finite. For a group $G$, we write $|G|$ for its order, $o(g)$ for the order of an element $g\in G$, $\langle X\rangle$ for the subgroup generated by a subset $X\subseteq G$, and $\mbox{exp}(G)$ for the exponent of $G$.
For a family $\mathcal{S}$ of subgroups of $G$, we also write
$\langle\mathcal{S}\rangle=\langle\bigcup_{H\in\mathcal{S}}H\rangle.$
The cyclic group of order $m$ is denoted by $C_m$. Let $U$ be a vector space over a field $F$. We write $\mbox{dim}_F(U)$ for the dimension of $U$, $\mbox{GL}_{F}(U)$ for the group of invertible $F$-linear transformations of $U$, $\mbox{End}_{F}(U)$ for the algebra of $F$-linear endomorphisms of $U$, and $\mbox{id}_U$ for the identity transformation on $U$. Moreover, for a set $Y\subseteq U$, we denote by $\mbox{span}_{F}(Y)$ the $F$-linear span of $Y$, that is, the set of all finite $F$-linear combinations of elements of $Y$. For a field $F$, $\mbox{char} F$ denotes its characteristic, and $F^\times$ denotes the multiplicative group of its nonzero elements. We write $\mathbb{F}_p$ for the finite field with $p$ elements, where $p$ is prime. The symbol $\cong$ denotes isomorphism between mathematical objects, when there is no doubt on which category they are considered in.

\section{Characterizations of consistent lattices}
In this section, we first characterize consistent lattices of finite length via their sublattices. We then characterize consistency for complete atomistic lattices of finite length using closure operators. 

Let $L$ be a lattice. For convenience, for any $x, y\in L$, if neither $x\leq y$ nor $y\leq x$ holds, then $x$ and $y$ are said to be incomparable, and we write $x\parallel y$. Otherwise, we use symbols $x\nparallel y$. Further background on lattice theory may be found in \cite{Crawley73,Birkhoff73,Davey,Stern99}.
\subsection{Consistency in lattices of finite length}
Although Stern \cite{Stern99} characterized consistency for semimodular lattices of finite length via their sublattices, analogous characterizations 
for arbitrary lattices of finite length have not yet been obtained. 
In what follows, we develop a theory characterizing the consistency 
of lattices of finite length in terms of their sublattices.

\begin{definition}\label{D01}
\emph{A subposet $P$ of a lattice $L$ is said to be operation-preserving if, for all $x, y\in P$,
$$x\wedge_P y \mbox{ exists }\Rightarrow x\wedge_P y=x\wedge_L y,$$
and
$$x\vee_P y \mbox{ exists }\Rightarrow x\vee_P y=x\vee_L y.$$}
\end{definition}

Every sublattice of $L$ is an operation-preserving subposet of $L$. However, the converse does not hold in general. 
For example, let $L$ denote the lattice of all subsets of the finite set $\{a, b, c\}$ ordered by set inclusion. 
Then the subposet $$P=(\{\{a\}, \{a,b\}, \{a,c\}, \{b,c\}, \{a,b,c\}\}, \subseteq)$$ is an operation-preserving 
subposet of $L$. However, $P$ is not a sublattice of $L$ since $\{a,b\}\wedge_L \{b,c\}=\{b\}\notin P$. 
For convenience, if $P$ is an operation-preserving subposet of $L$, then, whenever no confusion can arise, 
we write $\vee=\vee_L=\vee_P$ and $\wedge=\wedge_L=\wedge_P$, respectively.

\begin{lemma}\label{L01}
A lattice $L$ of finite length is consistent if and only if it does not contain an operation-preserving subposet of the form shown in Figure 1, with the two covers in $L$ indicated by double lines and the distinguished element $j\in J^{+}(L)$. 
\end{lemma}
\begin{center}
\begin{minipage}{11pc}
\setlength{\unitlength}{1pt}\begin{picture}(600,130)
\put(40,120){\circle{6}}\put(30,130){\makebox(0,0)[l]{\footnotesize$t\vee j$}}
\put(40,40){\circle{6}}\put(45,40){\makebox(0,0)[l]{\footnotesize$t$}}
\put(20,80){\circle{6}}\put(05,80){\makebox(0,0)[l]{\footnotesize$a$}}
\put(60,80){\circle{6}}\put(65,80){\makebox(0,0)[l]{\footnotesize$b$}}
\put(100,80){\circle{6}}\put(105,80){\makebox(0,0)[l]{\footnotesize$j \quad (j\in J^{+}(L))$}}
\put(38,42.3){\line(-1,2){17.5}}
\put(42,42.3){\line(1,2){17.5}}
\put(22,82.3){\line(1,2){17.5}}
\put(20,82.3){\line(1,2){18}}
\put(58,82.3){\line(-1,2){17.5}}
\put(60,82.3){\line(-1,2){18}}
\put(98,82){\line(-3,2){55}}
\put(38,23){$P$}
\put(-60,5){Figure 1. An operation-preserving subposet}
\end{picture}
\end{minipage}
\end{center}
\proof
Suppose that $L$ contains an operation-preserving subposet $P$ as illustrated in Figure 1. Then the element $t\vee j$ is not join-irreducible in the subinterval $[t,1]$. 
Hence, $j$ is inconsistent, which implies that $L$ is inconsistent.

Conversely, suppose that $L$ is inconsistent. Then there exist elements $j\in J^{+}(L)$ and $x\in L$ such 
that $x\vee j$ is not join-irreducible in the sublattice $[x,1]$. Hence there exist two distinct elements $u, v\in [x, x\vee j]$ such that 
\begin{equation}\label{eq+1}
 u, v<x\vee j\mbox{ and }u\vee v=x\vee j.   
\end{equation}
Since $L$ has finite length, we may choose $a$ and $b$ to be two maximal elements below $x\vee j$ 
such that $u\leq a\prec x\vee j$ and $v\leq b\prec x\vee j$.  By \eqref{eq+1}, $a\vee b=x\vee j$. 
We claim that $j\parallel a$ and $j\parallel b$. Indeed, if $j\leq a$ then $x\vee j\leq a<x\vee j$, a contradiction. 
If $a\leq j$, then $x\vee j\leq a\vee j=j\in J^{+}(L)$. Thus $x\leq j$, which yields that $x\vee j=j\in J([x,1])$, 
contrary to the fact that $x\vee j\notin J([x,1])$. Thus, $j\parallel a$. Similarly, $j\parallel b$. Put $a\wedge b=t$. 
Then $x\leq t$ since $x\leq a, b$. Hence
$$x\vee j\leq t\vee j=(a\wedge b)\vee j\leq a\vee j\leq x\vee j,$$
which implies that $x\vee j=t\vee j$. Therefore, the elements $j, a, b, t$ and $t\vee j$ form a subposet of $L$ of 
the form shown in Figure 1. Consequently, $L$ contains an operation-preserving subposet of the form shown in Figure 1. 
\endproof

\begin{example}\label{E01}
\emph{Consider the subgroup lattice $L(A_4)$ of the alternating group $A_4$, as depicted in Figure 2. It is straightforward to verify that $L(A_4)$ does not contain an operation-preserving subposet of the form shown in Figure 1. Therefore, by Lemma \ref{L01}, $L(A_4)$ is a consistent lattice. It is worth noting that $L(A_4)$ is neither an upper semimodular lattice nor a lower semimodular lattice.} 
\end{example}

\begin{center}
\begin{minipage}{11pc}
\setlength{\unitlength}{0.8pt}\begin{picture}(600,180)
\put(-80,80){\circle{6}}\put(-95,80){\makebox(0,0)[l]{\footnotesize$a$}}
\put(-40,80){\circle{6}}\put(-31,81){\makebox(0,0)[l]{\footnotesize$b$}}
\put(0,80){\circle{6}}\put(5,80){\makebox(0,0)[l]{\footnotesize$c$}}
\put(40,80){\circle{6}}\put(45,80){\makebox(0,0)[l]{\footnotesize$d$}}
\put(80,80){\circle{6}}\put(85,80){\makebox(0,0)[l]{\footnotesize$e$}}
\put(120,80){\circle{6}}\put(126,80){\makebox(0,0)[l]{\footnotesize$f$}}
\put(160,80){\circle{6}}\put(165,80){\makebox(0,0)[l]{\footnotesize$g$}}
\put(40,40){\circle{6}}\put(38,32){\makebox(0,0)[l]{\footnotesize$0$}}
\put(40,120){\circle{6}}\put(47,120){\makebox(0,0)[l]{\footnotesize$h$}}
\put(40,160){\circle{6}}\put(36,170){\makebox(0,0)[l]{\footnotesize$A_4$}}
\put(40,43){\line(0,1){34}}
\put(40,83){\line(0,1){34}}
\put(40,123){\line(0,1){34}}
\put(42,42){\line(1,1){36}}
\put(42,42){\line(2,1){74}}
\put(42,42){\line(3,1){115}}
\put(38,42){\line(-1,1){36}}
\put(38,42){\line(-2,1){74}}
\put(38,42){\line(-3,1){115}}
\put(2,82){\line(1,1){36}}
\put(78,82){\line(-1,1){36}}
\put(118,82){\line(-1,1){76}}
\put(158,82){\line(-3,2){115}}
\put(-38,82){\line(1,1){76}}
\put(-78,82){\line(3,2){116}}
\put(30,14){$L(A_4)$}
\put(-30,-5){Figure 2. The lattice $L(A_4)$}
\end{picture}
\end{minipage}
\end{center}

Let $L$ be a lattice of finite length. In the following, for any element $j\in J^+(L)$, we denote 
$$\Im(j)=\{t\in L: t\vee j\notin J([t,1])\}.$$ 
By the definition of consistency, if $j$ is inconsistent then $\Im(j)\neq \emptyset$, and hence $\Im(j)$ forms a subposet of $L$. 
Consequently, there exists at least one maximal element $t_0\in \Im(j)$ for any inconsistent element $j\in J^{+}(L)$ since $L$ is a lattice of finite length.

\begin{example}\label{E002}
\emph{Consider the finite lattice $L$ depicted in Figure 3. It is easy to verify that $j\in J^{+}(L)$ and $j$ is inconsistent. Moreover, we know that 
$$\Im(j)=\{t\in L: t\vee j\notin J([t,1])\}=\{t_0, t_1\}.$$
Therefore, $\Im(j)$ is a subposet of $L$ and $t_1$ is the maximal element in $\Im(j)$.}
\end{example}
\begin{center}
	\begin{minipage}{11pc}
\setlength{\unitlength}{0.8pt}\begin{picture}(600,210)
\put(-80,80){\circle{6}}\put(-95,80){\makebox(0,0)[l]{\footnotesize$j$}}
\put(-110,65){\makebox(0,0)[l]{\footnotesize$(j\in J^{+}(L))$}}
\put(0,120){\circle{6}}\put(5,120){\makebox(0,0)[l]{\footnotesize$a$}}
\put(40,80){\circle{6}}\put(45,80){\makebox(0,0)[l]{\footnotesize$t_0$}}
\put(80,120){\circle{6}}\put(85,120){\makebox(0,0)[l]{\footnotesize$t_1$}}
\put(80,160){\circle{6}}\put(86,160){\makebox(0,0)[l]{\footnotesize$b$}}
\put(120,160){\circle{6}}\put(125,160){\makebox(0,0)[l]{\footnotesize$c$}}
\put(40,40){\circle{6}}\put(38,32){\makebox(0,0)[l]{\footnotesize$0$}}
\put(40,200){\circle{6}}\put(38,210){\makebox(0,0)[l]{\footnotesize$1$}}
\put(40,43){\line(0,1){34}}
\put(80,123){\line(0,1){34}}
\put(38,42){\line(-3,1){115}}
\put(38,82){\line(-1,1){36}}
\put(42,82){\line(1,1){36}}
\put(82,122){\line(1,1){36}}
\put(78,162){\line(-1,1){36}}
\put(2,122){\line(1,2){38}}
\put(118,162){\line(-2,1){75}}
\put(-78,82){\line(1,1){116}}
\put(36,15){$L$}
\put(-65,3){Figure 3. A finite inconsistent lattice}
\end{picture}
\end{minipage}
\end{center}

\begin{lemma}\label{L02}
Let $L$ be a lattice of finite length, and let $j\in J^{+}(L)$ be inconsistent. 
If $t$ is a maximal element of $\Im(j)$, then $L$ contains an operation-preserving subposet 
of the form shown in Figure 1, with two covers indicated by double lines.
\end{lemma}
\proof 
Since $t\in \Im(j)$, we know that $t\vee j\notin J([t,1])$, which implies that $t\vee j\notin J([t,t\vee j])$. 
Thus, there exist at least two distinct elements $a, b\in [t, t\vee j]$ such that both $a\prec t\vee j$ and $b\prec t\vee j$ 
since $L$ is a lattice of finite length. 
Therefore, $a\wedge b\geq t, a\parallel b$ and $a\vee b=t\vee j$. We claim that $a\wedge b=t$. Otherwise, $t\vee j>a\wedge b> t$. 
Hence, $t\vee j\geq (a\wedge b)\vee j\geq t\vee j$, 
which implies that $$(a\wedge b)\vee j=t\vee j \notin J([a\wedge b,1]).$$ Therefore, $a\wedge b\in \Im(j)$, 
contrary to the fact that $t$ is a maximal element in $\Im(j)$. 
Consequently, the elements $a, b, j,t$ and $t\vee j$ form an operation-preserving subposet of the form shown in Figure 1.
\endproof

\begin{example}\label{E02}
\emph{Consider the finite lattice $L$ depicted in Figure 3 again. It is easy to verify that $t_1$ is maximal in $\Im(j)$ and $\{j,b, c, t_1, 1\}$ is 
an operation-preserving subposet of the form shown in Figure 1.}
\end{example}

\begin{remark}\label{R01}
\emph{The lattice $L$ depicted in Figure 3 also tells us that $t_0\in \Im(j)$ and $\{j,a, b, t_0, 1\}$ 
is also an operation-preserving subposet of the form shown in Figure 1. 
However, $t_0$ is not maximal in $\Im(j)$. Therefore, the converse of Lemma \ref{L02} does not hold in general.} 
\end{remark}

\begin{lemma}\label{T01}
A lattice of finite length is consistent if and only if it contains no sublattice of any of the three forms shown 
in Figure 4, with the distinguished element $j\in J^{+}(L)$.  
\end{lemma}
\begin{center}
\begin{minipage}{11pc}
\setlength{\unitlength}{1pt}\begin{picture}(600,180)
\put(140,40){\circle{6}}
\put(180,80){\circle{6}}
\put(100,80){\circle{6}}
\put(140,120){\circle{6}}
\put(180,120){\circle{6}}\put(185,120){\makebox(0,0)[l]{\footnotesize$j\quad (j\in J^{+}(L))$}}
\put(100,120){\circle{6}}
\put(140,160){\circle{6}}

\put(-20,40){\circle{6}}
\put(-60,80){\circle{6}}
\put(-100,120){\circle{6}}
\put(20,120){\circle{6}}\put(25,120){\makebox(0,0)[l]{\footnotesize$j\quad (j\in J^{+}(L))$}}
\put(-20,120){\circle{6}}
\put(-60,160){\circle{6}}

\put(-22,42){\line(-1,1){36}}
\put(-62,82){\line(-1,1){36}}
\put(-58,82){\line(1,1){36}}
\put(-18,42){\line(1,2){37.5}}
\put(-98,122){\line(1,1){36}}
\put(-22,122){\line(-1,1){36}}
\put(18,122){\line(-2,1){75}}

\put(138,42){\line(-1,1){36}}
\put(100,83){\line(0,1){34}}
\put(142,42){\line(1,1){36}}
\put(180,83){\line(0,1){34}}
\put(178,82){\line(-1,1){36}}
\put(102,82){\line(1,1){36}}
\put(102,122){\line(1,1){36}}
\put(178,122){\line(-1,1){36}}
\put(140,123){\line(0,1){34}}
\put(-23,25){$L_5$}
\put(135,25){$S_7$}
\end{picture}
\end{minipage}
\end{center}

\begin{center}
\begin{minipage}{11pc}
\setlength{\unitlength}{1pt}\begin{picture}(600,180)
\put(60,40){\circle{6}}
\put(100,80){\circle{6}}
\put(20,80){\circle{6}}
\put(60,120){\circle{6}}
\put(60,80){\circle{6}}
\put(100,120){\circle{6}}
\put(20,120){\circle{6}}
\put(60,160){\circle{6}}\put(67,160){\makebox(0,0)[l]{\footnotesize$j \quad (j\in J^{+}(L))$}}
\put(60,200){\circle{6}}
\put(58,42){\line(-1,1){36}}
\put(20,83){\line(0,1){34}}
\put(62,42){\line(1,1){36}}
\put(100,83){\line(0,1){34}}
\put(98,82){\line(-1,1){36}}
\put(22,82){\line(1,1){36}}
\put(22,122){\line(1,2){37.7}}
\put(98,122){\line(-1,2){37.7}}
\put(60,123){\line(0,1){34}}
\put(60,43){\line(0,1){34}}
\put(60,163){\line(0,1){34}}
\put(62,82){\line(1,1){36}}
\put(58,82){\line(-1,1){36}}
\put(57,25){$B_{3}^*$}
\put(-40,5){Figure 4. Three inconsistent sublattices}
\end{picture}
\end{minipage}
\end{center}
\proof
Suppose that $L$ is an inconsistent lattice of finite length. Then there exists an inconsistent element $j\in J^{+}(L)$. 
Hence, there exists a maximal element $t_0\in \Im(j)$. 
By Lemma \ref{L02}, $L$ contains an operation-preserving subposet $P=\{a, b, t_0, j, t_0\vee j\}$ of the form depicted in Figure 1. 
Therefore, the following equations hold:
\begin{equation}\label{Eq01}
a\parallel b, a\prec t_0\vee j, b\prec t_0\vee j, a\wedge b=t_0 \mbox{ and } a\vee b=t_0\vee j=a\vee j=b\vee j.
\end{equation}
The proof proceeds in two steps as follows.

A. Suppose that $a\wedge j \nparallel b\wedge j$. Then since $a$ and $b$ play symmetric roles, there are two cases.

Case 1. Suppose that $a\wedge j=b\wedge j$. Then $a\wedge j=b\wedge j=t_0\wedge j$. 
Thus, from \eqref{Eq01}, we know that $L$ contains a sublattice of the form $L_5$ shown in Figure 4.

Case 2. Suppose that $a\wedge j<b\wedge j$. Then $a\wedge j=a\wedge b\wedge j=t_0\wedge j$. 
We first claim that $b\wedge j\parallel a$. Indeed, if $a<b\wedge j$, then $a<b\wedge j\leq b$, contrary to the fact that $a\parallel b$. 
If $b\wedge j\leq a$, then $b\wedge j\leq b\wedge a=t_0$, and hence $b\wedge j\leq t_0\wedge j=a\wedge j< b\wedge j$, a contradiction. Therefore, $$b\wedge j\parallel a.$$ 
Note that $a\prec t_0\vee j$. Then 
\begin{equation}\label{eq+2}
a\vee (b\wedge j)=t_0\vee j.    
\end{equation} 
Put $b_1=(b\wedge j)\vee t_0$. Then $t_0\leq b_1\leq b$. Moreover, $b_1\nleq a$, for otherwise $b\wedge j\leq b_1\leq a$, 
contrary to $b\wedge j\parallel a$. Therefore, two subcases can be distinguished.

Subcase 1. If $b_1=b$, then the set $\{a, b, j, b\wedge j, t_0, t_0\vee j, t_0\wedge j\}$ forms a sublattice of the type $S_7$ shown in Figure 4 by \eqref{Eq01} and \eqref{eq+2}.

Subcase 2. If $b_1<b$ then $t_0\leq b_1< b$. 
Moreover, we claim that $t_0<b_1$. Otherwise $b_1=t_0$, and hence $b\wedge j\le t_0\wedge j=a\wedge j$, a contradiction.
Thus, by \eqref{Eq01} and \eqref{eq+2}, we have that $$a\wedge b_1=a\wedge b=t_0,\ a\vee b_1=t_0\vee j,\ b_1\vee j=b\vee j=t_0\vee j \mbox{ and }b_1\wedge j=b\wedge j.$$ 
Therefore, $L$ contains a sublattice $U$ as shown in Figure 5.
\begin{center}
\begin{minipage}{11pc}
\setlength{\unitlength}{0.8pt}\begin{picture}(600,180)
\put(60,40){\circle{6}}\put(52,30){\makebox(0,0)[l]{\footnotesize$t_0\wedge j$}}
\put(100,80){\circle{6}} \put(105,80){\makebox(0,0)[l]{\footnotesize$b\wedge j$}}
\put(20,80){\circle{6}}\put(6,80){\makebox(0,0)[l]{\footnotesize$t_0$}}
\put(60,120){\circle{6}}\put(46,122){\makebox(0,0)[l]{\footnotesize$b_1$}}
\put(100,120){\circle{6}}\put(105,120){\makebox(0,0)[l]{\footnotesize$j \quad (j\in J^{+}(L))$}}
\put(20,120){\circle{6}}\put(8,120){\makebox(0,0)[l]{\footnotesize$a$}}
\put(60,160){\circle{6}}\put(48,170){\makebox(0,0)[l]{\footnotesize$t_0\vee j$}}
\put(60,140){\circle{6}}\put(65,140){\makebox(0,0)[l]{\footnotesize$b$}}
\put(58,42){\line(-1,1){36}}
\put(20,83){\line(0,1){34}}
\put(62,42){\line(1,1){36}}
\put(100,83){\line(0,1){34}}
\put(98,82){\line(-1,1){36}}
\put(22,82){\line(1,1){36}}
\put(22,122){\line(1,1){36}}
\put(98,122){\line(-1,1){36}}
\put(60,123){\line(0,1){14}}
\put(60,143){\line(0,1){14}}
\put(-20,5){Figure 5. A sublattice $U$ of $L$}
\end{picture}
\end{minipage}
\end{center}
It is clear that $U$ contains a sublattice of the form $S_7$ shown in Figure 4. 
Consequently, $L$ contains a sublattice of the form $S_7$ shown in Figure 4..

B. Suppose that $a\wedge j\parallel b\wedge j$. Then set $(a\wedge j)\vee (b\wedge j)=r$. Thus $r\leq j$. 
Note that $a\wedge j<j$ and $b\wedge j<j$. Hence, it follows that $r<j$ since $j\in J^{+}(L)$. We claim that
$$r\parallel a \mbox{ and } r\parallel b.$$
Indeed, if $r\leq a$ then $b\wedge j\leq r \leq a$, which together with $b\wedge j\leq r<j$ implies that $b\wedge j\leq a\wedge j$, contrary to the assumption that $a\wedge j\parallel b\wedge j$.
If $a<r$ then $a<r\leq t_0\vee j$, which together with $a\prec t_0\vee j$ implies that $r=t_0\vee j>j>r$, a contradiction. 
Thus, $r\parallel a$. Since $a$ and $b$ play symmetric roles, we similarly conclude that $r\parallel b$. Therefore,
$$(b\wedge j)\vee a=(b\wedge j)\vee (a\wedge j) \vee a=r\vee a=t_0\vee j$$ since $r\parallel a$ and $a\prec t_0\vee j$. Similarly,
$$(a\wedge j)\vee b=r\vee b=t_0\vee j.$$ The proof proceeds by considering the following four cases.

Case a. Suppose that $(a\wedge j)\vee t_0=a$ and $(b\wedge j)\vee t_0=b$. Then the elements $t_0\wedge j, a\wedge j, b\wedge j, t_0, r, j, a, b$ and $ t_0\vee j$ constitute a sublattice of the form $B_{3}^{*}$ shown in Figure 4.  

Case b. Suppose that $(a\wedge j)\vee t_0=a$ and $(b\wedge j)\vee t_0<b$. Put $[(b\wedge j)\vee t_0]\vee r=k$.  Then $k\leq t_0\vee j$. 
We first claim that $k\parallel b$ when $k<t_0\vee j$. Indeed, if $b<k$ then $b<k<t_0\vee j$, 
contrary to the fact that $b\prec t_0\vee j$. If $k\leq b$ then $r\leq k\leq b$, contrary to the fact that $b\parallel r$. Therefore, $k\parallel b$ when $k<t_0\vee j$.
 We now prove that $k=t_0\vee j$. Assume, to the contrary, that $k<t_0\vee j$. Observe that $k\vee b=t_0\vee j$ since $k\parallel b$ and $b\prec t_0\vee j$. 
Moreover, $\{k,b\}\subseteq [t_0\vee(b\wedge j),t_0\vee j]\subseteq [t_0\vee(b\wedge j),1]$. 
Thus, $$[t_0\vee(b\wedge j)]\vee j=t_0\vee j\notin J([t_0\vee(b\wedge j),1]).$$
By the definition of $\Im(j)$, we have $t_0\vee (b\wedge j)\in \Im(j)$. Therefore, $t_0\vee (b\wedge j)= t_0$ since $t_0$ is the maximal element of $\Im(j)$ and $t_0\leq t_0\vee (b\wedge j)$. 
Hence $b\wedge j\leq t_0$. This implies that $b\wedge j=b\wedge j\wedge t_0\leq t_0\wedge j=a\wedge b\wedge j\leq a\wedge j$, contrary to the assumption that $b\wedge j\parallel a\wedge j$. Consequently, 
\begin{equation}\label{eq+3}
[(b\wedge j)\vee t_0]\vee r=k = t_0 \vee j.   
\end{equation}
It follows from the above that $L$ contains a sublattice $R$ shown in Figure 6.
\begin{center}
\begin{minipage}{11pc}
\setlength{\unitlength}{0.8pt}\begin{picture}(600,180)
\put(60,40){\circle{6}}\put(52,30){\makebox(0,0)[l]{\footnotesize$t_0\wedge j$}}
\put(100,80){\circle{6}} \put(105,80){\makebox(0,0)[l]{\footnotesize$b\wedge j$}}
\put(20,80){\circle{6}}\put(-10,80){\makebox(0,0)[l]{\footnotesize$a\wedge j$}}
\put(60,120){\circle{6}}\put(48,122){\makebox(0,0)[l]{\footnotesize$r$}}
\put(60,80){\circle{6}}\put(44,82){\makebox(0,0)[l]{\footnotesize$t_0$}}
\put(100,120){\circle{6}}\put(105,120){\makebox(0,0)[l]{\footnotesize$t_0\vee (b\wedge j) \quad (j\in J^{+}(L))$}}
\put(80,140){\circle{6}}\put(85,140){\makebox(0,0)[l]{\footnotesize$b$}}
\put(20,120){\circle{6}}\put(8,120){\makebox(0,0)[l]{\footnotesize$a$}}
\put(60,160){\circle{6}}\put(48,170){\makebox(0,0)[l]{\footnotesize$t_0\vee j$}}
\put(60,140){\circle{6}}\put(65,140){\makebox(0,0)[l]{\footnotesize$j$}}
\put(58,42){\line(-1,1){36}}
\put(20,83){\line(0,1){34}}
\put(62,42){\line(1,1){36}}
\put(100,83){\line(0,1){34}}
\put(98,82){\line(-1,1){36}}
\put(22,82){\line(1,1){36}}
\put(22,122){\line(1,1){36}}
\put(98,122){\line(-1,1){16}}
\put(78,142){\line(-1,1){16}}
\put(60,123){\line(0,1){14}}
\put(60,43){\line(0,1){34}}
\put(62,82){\line(1,1){36}}
\put(58,82){\line(-1,1){36}}
\put(60,143){\line(0,1){14}}
\put(10,5){Figure 6. A sublattice $R$}
\end{picture}
\end{minipage}
\end{center}
As $R$ is a sublattice of $L$, the set $\{a, t_0\vee(b\wedge j) , j, b\wedge j, r, a\wedge j, t_0, t_0\vee j, t_0\wedge j\}$ constitutes a sublattice of the form $B_{3}^{*}$.

Case c. Suppose that $(a\wedge j)\vee t_0<a$ and $(b\wedge j)\vee t_0=b$. Then since $a$ and $b$ play symmetric roles, similarly to Case b, $L$ must also contain a sublattice of the form $B_{3}^{*}$.

Case d. Suppose that $(a\wedge j)\vee t_0<a$ and $(b\wedge j)\vee t_0<b$. Then similarly to the proof of \eqref{eq+3}, we obtain that $$(a\wedge j)\vee t_0\vee (b\wedge j)=[(a\wedge j)\vee t_0] \vee r=[(b\wedge j)\vee t_0]\vee r=t_0\vee j.$$ Therefore, we have the following sublattice $F$ of $L$ as shown in Figure 7.
\begin{center}
\begin{minipage}{11pc}
\setlength{\unitlength}{0.8pt}\begin{picture}(600,180)
\put(60,40){\circle{6}}\put(52,30){\makebox(0,0)[l]{\footnotesize$t_0\wedge j$}}
\put(100,80){\circle{6}} \put(105,80){\makebox(0,0)[l]{\footnotesize$b\wedge j$}}
\put(20,80){\circle{6}}\put(-10,80){\makebox(0,0)[l]{\footnotesize$a\wedge j$}}
\put(60,120){\circle{6}}\put(48,122){\makebox(0,0)[l]{\footnotesize$r$}}
\put(60,80){\circle{6}}\put(44,82){\makebox(0,0)[l]{\footnotesize$t_0$}}
\put(100,120){\circle{6}}\put(105,120){\makebox(0,0)[l]{\footnotesize$t_0\vee (b\wedge j) \quad (j\in J^{+}(L))$}}
\put(80,140){\circle{6}}\put(85,140){\makebox(0,0)[l]{\footnotesize$b$}}
\put(40,140){\circle{6}}\put(29,140){\makebox(0,0)[l]{\footnotesize$a$}}
\put(20,120){\circle{6}}\put(-44,120){\makebox(0,0)[l]{\footnotesize$t_0\vee (a\wedge j)$}}
\put(60,160){\circle{6}}\put(48,170){\makebox(0,0)[l]{\footnotesize$t_0\vee j$}}
\put(60,140){\circle{6}}\put(65,140){\makebox(0,0)[l]{\footnotesize$j$}}
\put(58,42){\line(-1,1){36}}
\put(20,83){\line(0,1){34}}
\put(62,42){\line(1,1){36}}
\put(100,83){\line(0,1){34}}
\put(98,82){\line(-1,1){36}}
\put(22,82){\line(1,1){36}}
\put(22,122){\line(1,1){16}}
\put(42,142){\line(1,1){16}}
\put(98,122){\line(-1,1){16}}
\put(78,142){\line(-1,1){16}}
\put(60,123){\line(0,1){14}}
\put(60,43){\line(0,1){34}}
\put(62,82){\line(1,1){36}}
\put(58,82){\line(-1,1){36}}
\put(60,143){\line(0,1){14}}
\put(10,5){Figure 7. A sublattice $F$}
\end{picture}
\end{minipage}
\end{center}
Thus, the set $\{t_0\vee(a\wedge j), t_0\vee(b\wedge j), j, b\wedge j, r, a\wedge j, t_0, t_0\vee j, t_0\wedge j\}$ constitutes a sublattice of the form $B_{3}^{*}$.

Combining A and B, every inconsistent lattice of finite length contains a sublattice of one of the forms $L_5$, $S_7$, or $B_{3}^{*}$ in Figure 4.

Conversely, suppose that $L$ contains a sublattice of the form $L_5$, $S_7$, or $B_{3}^{*}$ shown in Figure 4. 
In each of these cases, there exists an element $t\in L$ such that $t\vee j\notin J([t,1])$. Hence, the nonzero join-irreducible element $j$ is inconsistent in $L$, 
and therefore $L$ is inconsistent.
\endproof
\quad\\
\textbf{Proof of Theorem \ref{TT01}:} Theorem \ref{TT01} follows from Lemmas \ref{L01} and \ref{T01}.\endproof 
\begin{corollary}\label{C02}
A finite planar lattice is consistent if and only if it contains no sublattice of either of the two forms $L_5$ and $S_7$ shown in Figure 4, with the distinguished element $j\in J^{+}(L)$.
\end{corollary}
\proof
From Theorem 1 in Kelly and Rival \cite{Kelly}, a finite planar lattice $L$ contains no sublattice isomorphic to the nonplanar lattice $B_{3}^{*}$. Therefore, Corollary \ref{C02} holds by Lemma \ref{T01}.
\endproof

\subsection{Consistency in complete atomistic lattices}
Let $L$ be a finite atomistic lattice and let $\mbox{cl}$ be the canonical closure operator associated with $L$. 
Then, $L\cong \mathcal{C}(A(L), \mbox{cl})$ is geometric 
if $\mbox{cl}$ satisfies EP \cite{Crawley73, Birkhoff73}. 
Every geometric lattice is consistent; however, the converse does not hold \cite{Stern99}.
We next establish equivalent characterizations of consistency for atomistic lattices 
via closure operators satisfying CEP.

\begin{remark}\label{r00}
\emph{If $L$ is a complete atomistic lattice of finite length, then $J^{+}(L)=A(L)$. Thus, by Theorem \ref{TT01}, a complete atomistic lattice of finite length is consistent if and only if it contains no sublattice of the form $L_{5}$ shown in Figure 4, with the distinguished element $j\in J^{+}(L)$.}    
\end{remark}

\begin{corollary}\label{C01}
A complete atomistic lattice $L$ of finite length is consistent if and only if it contains no $0$-sublattice 
of the form $L_{5}^{u}$ shown in Figure 8 with three covers in $L$ indicated by double lines and the distinguished element $j\in J^{+}(L)$.      
\end{corollary}
\begin{center}
\begin{minipage}{11pc}
\setlength{\unitlength}{0.8pt}\begin{picture}(600,180)
\put(60,40){\circle{6}}\put(58,31){\makebox(0,0)[l]{\footnotesize$0$}}
\put(20,80){\circle{6}}\put(8,80){\makebox(0,0)[l]{\footnotesize$t_0$}}
\put(-20,120){\circle{6}}
%\put(-30,120){\makebox(0,0)[l]{\footnotesize$a$}}
\put(100,120){\circle{6}}\put(105,120){\makebox(0,0)[l]{\footnotesize$j \quad (j\in J^{+}(L)=A(L))$}}
\put(60,120){\circle{6}}
%\put(65,120){\makebox(0,0)[l]{\footnotesize$b$}}
\put(20,160){\circle{6}}\put(8,170){\makebox(0,0)[l]{\footnotesize$t_0\vee j$}}
\put(58,42){\line(-1,1){36}}
\put(60,122.5){\line(-1,1){36.5}}
\put(18,82){\line(-1,1){36}}
\put(22,82){\line(1,1){36}}
\put(-20,122.5){\line(1,1){36.5}}
\put(62,42){\line(1,2){37.5}}
\put(60,42.5){\line(1,2){38}}
\put(-18,122){\line(1,1){36}}
\put(58,122){\line(-1,1){36}}
\put(98,122){\line(-2,1){75}}
\put(-23,5){Figure 8. A $0$-sublattice $L_{5}^{u}$}
\end{picture}
\end{minipage}
\end{center}
\proof
Suppose that $L$ is inconsistent. Since $L$ is an atomistic lattice of finite length, $J^{+}(L)=A(L)$. By Theorem \ref{TT01}, it follows that $L$ contains a $0$-sublattice of the form $L_{5}^{u}$. By Remark \ref{r00}.

The converse follows immediately.
\endproof

\begin{lemma}\label{L04}
A complete atomistic lattice $L$ of finite length is consistent if and only if it contains no $0$-sublattice of the form $L_{5}^{l}$ shown in Figure 9 with three covers in $L$ indicated by double lines and the distinguished element $j\in J^{+}(L)$.
\end{lemma}
\begin{center}
\begin{minipage}{11pc}
\setlength{\unitlength}{0.8pt}\begin{picture}(600,170)
\put(60,40){\circle{6}}\put(58,31){\makebox(0,0)[l]{\footnotesize$0$}}
\put(20,80){\circle{6}}\put(8,80){\makebox(0,0)[l]{\footnotesize$t_0$}}
\put(-20,120){\circle{6}}
%\put(-30,120){\makebox(0,0)[l]{\footnotesize$a$}}
\put(100,120){\circle{6}}\put(105,120){\makebox(0,0)[l]{\footnotesize$j \quad (j\in J^{+}(L)=A(L))$}}
\put(60,120){\circle{6}}
%\put(65,120){\makebox(0,0)[l]{\footnotesize$b$}}
\put(20,160){\circle{6}}\put(8,170){\makebox(0,0)[l]{\footnotesize$t_0\vee j$}}
\put(58,42){\line(-1,1){36}}
\put(20,82.5){\line(-1,1){36.5}}
\put(18,82){\line(-1,1){36}}
\put(22,82){\line(1,1){36}}
\put(20,82.5){\line(1,1){36.5}}
\put(62,42){\line(1,2){37.5}}
\put(60,42.5){\line(1,2){38}}
\put(-18,122){\line(1,1){36}}
\put(58,122){\line(-1,1){36}}
\put(98,122){\line(-2,1){75}}
\put(-23,5){Figure 9. A $0$-sublattice $L_{5}^{l}$}
\end{picture}
\end{minipage}
\end{center}
\proof
By Remark \ref{r00}, if $L$ contains a $0$-sublattice of the form $L_{5}^{l}$, then $L$ is inconsistent.

 Conversely, suppose that $L$ is inconsistent. 
Then, by Corollary \ref{C01}, $L$ must contain a $0$-sublattice of the form $L_{5}^{u}$. Thus, there exist six elements $a, b, t_0, j, t_0\vee j$ and $0$ in $L$ such that 
\begin{equation}\label{eq+4}
a\vee b=t_0\vee j, t_0=a\wedge b, t_0\vee j\succ a>t_0,t_0\vee j\succ b>t_0  
\end{equation}
where $j\in J^{+}(L)$ and $t_0$ is maximal in $\Im(j)$.
It is evident that if $a\succ t_0$ and $b\succ t_0$, then $L$ contains a $0$-sublattice of the form $L_{5}^{l}$. It therefore remains to consider the following three cases.

Case (i). If $a\succ t_0$ and $b\nsucc t_0$, then there exists an element $c$ satisfying $b>c\succ t_0$. Therefore, by \eqref{eq+4}, the elements $0, t_0, a, c, t_0\vee j$ and $j$ constitute a sublattice of the form $L_{5}^{l}$ in Figure 9.

Case (ii). If $a\nsucc t_0$ and $b\succ t_0$, the argument proceeds analogously to Case (i).

Case (iii). If $a\nsucc t_0$ and $b\nsucc t_0$, then there exist two elements $c,d\in L$ such that $b>c\succ t_0$ and $a>d\succ t_0$. We claim that $$c\vee d=t_0\vee j.$$ 
Indeed, if $c\vee d<t_0\vee j$ then $b\nless c\vee d$ since $b\prec t_0\vee j$. Observe that $c\vee d\nleq b$, since otherwise $t_0\prec d\leq (c\vee d)\wedge a\leq b\wedge a=t_0$, a contradiction. 
Thus $(c\vee d)\parallel b$, which together with $t_0\vee j\succ b$ implies that $(c\vee d)\vee b=t_0\vee j$. Note that $c\vee j=t_0\vee j$. Hence, $c\in \Im(j)$ by the definition of $\Im(j)$. 
Therefore, $c=t_0$ because $t_0$ is maximal in $\Im(j)$, contrary to the fact that $t_0\prec c$. 
Therefore, $c\vee d=t_0\vee j$. Consequently, the elements $c, d, t_0, t_0\vee j, j$ and $0$ constitute a $0$-sublattice of the form $L^{l}_{5}$ shown in Figure 9.
\endproof

Below, we discuss the consistency of a complete atomistic lattice $L$ of finite length via the closure operator on $A(L)=J^{+}(L)$. 

\begin{proposition}\label{P01}
If a closure operator $\emph{cl}$ on a set $I$ satisfies EP, then it also satisfies CEP.
\end{proposition}
\proof
Suppose that $\mbox{cl}$ satisfies EP but not CEP. Then there exist $X\subseteq I$ and $x, y, z\in I$ such that
$$y\in \mbox{cl}(X\cup\{x,z\}) \setminus [\mbox{cl}(X\cup\{x\})\cup \mbox{cl}(X\cup \{z\})]$$  
and 
$$\mbox{cl}(X\cup\{x,z\})=\mbox{cl}(X\cup\{y\}).$$
Hence, $y\in \mbox{cl}(X\cup\{x,z\})\setminus \mbox{cl}(X\cup\{x\})$ and $y\in \mbox{cl}(X\cup\{x,z\})\setminus \mbox{cl}(X\cup\{z\})$. By EP, 
\begin{equation}\label{eq+5}
z\in \mbox{cl}(X\cup\{x,y\})\setminus \mbox{cl}(X\cup\{x\}) \mbox{ and } x\in \mbox{cl}(X\cup\{y,z\})\setminus \mbox{cl}(X\cup\{z\}). \end{equation}
Because $\mbox{cl}(X\cup\{x,z\})= \mbox{cl}(X\cup\{y\})$, we know that $$\mbox{cl}(X\cup\{x,y\})=\mbox{cl}(X\cup\{z,y\})=\mbox{cl}(X\cup\{y\}).$$
Hence, by \eqref{eq+5}, it follows that $x, z\in \mbox{cl}(X\cup \{y\}) \setminus \mbox{cl}(X)$ since $\mbox{cl}(X) \subseteq \mbox{cl}(X\cup\{x\})$ and $\mbox{cl}(X) \subseteq \mbox{cl}(X\cup\{z\})$. By EP again,
$$y\in \mbox{cl}(X\cup \{x\}) \setminus \mbox{cl}(X)\mbox{ and } y\in \mbox{cl}(X\cup \{z\}) \setminus \mbox{cl}(X),$$
contrary to the fact that $y\notin \mbox{cl}(X\cup \{x\})\cup \mbox{cl}(X\cup \{z\})$.
Consequently, if a closure operator $\mbox{cl}$ on a set $I$ satisfies EP then it also satisfies CEP.
\endproof

\quad\\
\textbf{Proof of Theorem \ref{T04}:} By Corollary \ref{C01} and Lemma \ref{L04}, (iii) $\Longleftrightarrow$ (i) $\Longleftrightarrow$ (iv). It therefore remains only to prove (iii) $\Longrightarrow$ (ii) and (ii) $\Longrightarrow$ (iv).

(iii) $\Longrightarrow$ (ii). Suppose that there exist $X\subseteq A(L)$ and $x, y, z\in A(L)$ such that
\begin{equation}\label{HE1}
y\in \mbox{cl}(X\cup\{x,z\}) \setminus [\mbox{cl}(X\cup\{x\})\cup \mbox{cl}(X\cup \{z\})]    
\end{equation} 
and 
\begin{equation}\label{HE2}
\mbox{cl}(X\cup\{x,z\})=\mbox{cl}(X\cup\{y\}).
\end{equation}
We first claim that $X\neq \emptyset$. Otherwise, $\mbox{cl}(\{x,z\})=\mbox{cl}(\{y\})$. Then by \eqref{eqqq-1} and 
$A(L)=J^{+}(L)$, $$\{x,z\}\subseteq A(x\vee z)=A(y)=\{y\}.$$
Thus $x=y=z$, which implies that $$\mbox{cl}(X\cup\{x,z\}) \setminus [\mbox{cl}(X\cup\{x\})\cup \mbox{cl}(X\cup \{z\})]=\emptyset,$$ a contradiction.
Set $t=\bigvee X$. Then, by \eqref{HE2}, $t\vee y=t\vee x\vee z$. Let $p=t\vee y$. Since $y\notin \mbox{cl}(X\cup\{x\})$, it follows that $t\vee x<p=t\vee x\vee z$. Similarly, we obtain $t\vee z<p$. Thus, there exist $a, b\in L$ such that $t\vee x\leq a\prec p$ and $t\vee z\leq b\prec p$. If $a=b$, then $p=t\vee x\vee z=a\vee b=a=b\prec p$, a contradiction. Therefore, $a\parallel b$ and $a\vee b=p$.

Let $t_0=a\wedge b$. Then $t\leq t_0< p$ since $t\leq a\prec p$ and $t\leq b\prec p$. We claim that $y\nleq a$. Otherwise, $p=t\vee y\leq a\prec p$, a contradiction. Similarly, $y\nleq b$. Thus, $y\nleq a\wedge b=t_0$. Hence, $y\wedge t_0=0$ and $p=t\vee y\leq t_0\vee y\leq p$ since $y\in A(L)$ and $t\leq t_0<p$. Then it follows that $t_0\vee y=p$. 
Therefore, $$a\wedge b=t_0, a\vee b=p=t_0\vee y, t_0\wedge y=0, a\prec p \mbox{ and } b\prec p.$$ 
Note that $t>0$ and $y\wedge b=y\wedge a=0$ since $X\neq \emptyset$ and $y\nleq a, b$. Put $y=j$. Thus, the elements $0, t_0, a, b, y$ and $t_0\vee y$ form a $0$-sublattice of the form $L_{5}^{u}$ shown in Figure 8 where $y=j\in A(L)=J^{+}(L)$. This contradicts (iii). Therefore the canonical closure operator $\mbox{cl}$ satisfies CEP.

(ii) $\Longrightarrow$ (iv). Suppose that $L$ contains a $0$-sublattice of the form $L_{5}^{l}$ shown in Figure 9. Then there exist $j=y\in A(L)=J^{+}(L)$ and $t_0, a, b, t_0\vee y, 0\in L$ such that $$a\wedge b=t_0, a\vee b=t_0\vee y, a\wedge y=b\wedge y=t_0\wedge y=0, a\parallel b \mbox{ and } t_0\prec a, b.$$
Set $X=\{q\in A(L): q\leq t_0\}$ and $p=a\vee b=t_0\vee y$. 
Since $L$ is a complete atomistic lattice of finite length, there exists $x\in A(L)$ such that $x\leq a$ and $x\nleq t_0$, and hence $a=t_0\vee x$ because $t_0\prec a$. Similarly, there exists $z\in A(L)$ such that $z\leq b$, $z\nleq t_0$ and $b=t_0\vee z$. Thus $p=a\vee b=t_0\vee z\vee x=t_0\vee y$, 
which together with \eqref{eqqq-1} and $A(L)=J^{+}(L)$ implies that $$\mbox{cl}(X\cup\{x,z\})=\mbox{cl}(X\cup\{y\}).$$

We claim that $y\notin \mbox{cl}(X\cup \{x\})$. Otherwise, by \eqref{eqqq-1}, $y\leq t_0\vee x=a$, which implies $p=t_0\vee y\leq a<p$, a contradiction. Similarly, $y\notin \mbox{cl}(X\cup \{z\})$. Note that $y\in \mbox{cl}(X\cup \{x,z\})$ since $\mbox{cl}(X\cup\{x,z\})=\mbox{cl}(X\cup\{y\})$. Therefore, there exist $X\subseteq A(L)$ and $x, y, z\in A(L)$ such that
$$y\in \mbox{cl}(X\cup\{x,z\}) \setminus [\mbox{cl}(X\cup\{x\})\cup \mbox{cl}(X\cup \{z\})]$$  and  $$\mbox{cl}(X\cup\{x,z\})=\mbox{cl}(X\cup\{y\}),$$ a contradiction. 

This completes the proof of Theorem \ref{T04}. 
\endproof

\begin{remark}\label{R+11}
\emph{The converse of Proposition \ref{P01} fails. Indeed, the lattice $L(A_4)$ is a finite consistent atomistic lattice. 
Hence, by Theorem \ref{T04}, its canonical closure operator $\mbox{cl}$ satisfies CEP. However, we claim that $\mbox{cl}$ does not satisfy EP. Otherwise, 
$L(A_4)\cong \mathcal{C}(A(L(A_4)), \mbox{cl})$ is geometric, contrary to the fact that $L(A_4)$ is not semimodular. Therefore, by Proposition \ref{P01}, CEP is strictly weaker than EP.} 
\end{remark}

We now characterize the behavior of consistency under direct products, which will be used in the next section.
\begin{lemma}\label{Lem4.0}
Let $L_1$ and $L_2$ be two lattices of finite length. Then $L_1\times L_2$ is consistent if and only if both $L_1$ and $L_2$ are consistent.
\end{lemma}
\proof
From \cite{Stern99}, it follows that $L_1\times L_2$ is consistent whenever both $L_1$ and $L_2$ are consistent. 

Conversely, suppose that $L_1\times L_2$ is consistent. Write $0_i$ and $1_i$ for the least and greatest elements of $L_i$, respectively.
If $j\in J(L_1)$ and $x\in L_1$, then $(j,0_2)$ is join-irreducible in $L_1\times L_2$. Hence $$(x,0_2)\vee (j,0_2)=(x\vee j,0_2)\in J([(x,0_2), (1_1,1_2)]).$$ Therefore, $x\vee j\in J([x,1_1])$, which implies that $L_1$ is consistent. Similarly, we obtain that $L_2$ is also consistent.  \endproof

\section{Structure of finite groups}
As shown in \cite{Rotman} and \cite{Stern99}, $A_4\cong C^2_{2} \rtimes_{\theta} C_3$, and $L(A_4)$ is both consistent and atomistic. This example motivates the study of finite semidirect products whose subgroup lattices are consistent. In this section, we first characterize the atomistic property of the subgroup lattice $L(G)$ in terms of $\mbox{exp}(G)$, thereby establishing the connection needed to apply Theorem \ref{T04}. We then investigate $G=V\rtimes_\theta H$ to determine how consistency constrains the acting group $H$ and the associated $\mathbb{F}_pH$-module $V$. Our main results provide a complete characterization for cyclic actions and were extended to faithful actions of finite groups.

\subsection{Atomistic subgroup lattices}
A subgroup $M$ of $G$ is maximal if $M$ is a proper subgroup of $G$ and there is no subgroup $H$ such that $M<H<G$. The following proposition follows from the standard fact that a nontrivial finite group has a unique maximal subgroup if and only if it is a cyclic $p$-group for some prime $p$ \cite{Rotman}.
\begin{proposition}\label{P02}
Let $G$ be a finite group. Then $$J^{+}(L(G))=\{H\leq G: H \mbox{ is a nontrivial cyclic } p\mbox{-group}\}.$$  
Consequently, $$A(L(G))=\{H\leq G: |H| \mbox{ is prime}\}.$$
\end{proposition}
\proof
Suppose that $H\in J^{+}(L(G))$. Then there exists a unique 
lower cover $M$ of $H$ in $L(G)$. Thus $M$ is the unique maximal 
subgroup of $H$. It follows that $H\cong C_{p^n}$ for some prime $p$, i.e., $$J^{+}(L(G))\subseteq\{H\leq G: H \mbox{ is a nontrivial cyclic } p\mbox{-group}\}. $$
Conversely, if $H$ is a nontrivial cyclic $p$-group, then $H$ has 
a unique maximal subgroup and is therefore join-irreducible. 

Furthermore, as atoms of $L(G)$ are precisely the minimal nontrivial 
subgroups and $A(L(G))\subseteq J^{+}(L(G))$, 
these atoms are exactly the subgroups of prime order.
\endproof

\begin{proposition}\label{P03}
Let $G$ be a finite group. Then $L(G)$ is atomistic if and only if $\emph{exp}(G)$ is square-free.   
\end{proposition}
\proof
Suppose first that $L(G)$ is atomistic and $\mbox{exp}(G)$ is not square-free. 
Then there exists a prime $p$ such that $p^2\mid \mbox{exp}(G)$. 
Hence, there exists an element $g\in G$ such that $p^2\mid o(g)$. 
Let $o(g)=n$ and $K=\langle g^{n/p^2}\rangle$.
Then $|K|=p^2$, and therefore $K\cong C_{p^2}$. 
By Proposition \ref{P02}, $K\in J^{+}(L(G))$. 
Hence, the subgroup $K$ has a unique 
subgroup $K_1$ of order $p$. This implies that  
$$\{1_G\}\prec K_1\prec K.$$
Therefore, $\bigvee \{O\in A(L(G)): O\leq K\}=K_1<K$, a contradiction.

Conversely, suppose that $\mbox{exp}(G)$ is square-free. Let $H$ be a subgroup of $G$. Then, for any $h\in H$, $o(h)$ is square-free since $o(h)\mid \mbox{exp}(G)$. Note that there is nothing to prove when $h=1_G$. Assume $h\neq 1_G$. Thus, there exist $r$ distinct primes $p_1, p_2,\cdots,p_r$ such that $o(h)=p_1p_2\cdots p_r$. For each $i\in \{1,2,\cdots,r\}$, let $a_i=o(h)/p_i$. Then
$$d=\mbox{gcd}\{a_1,a_2,\cdots,a_r\}=1.$$
Hence, there exist $r$ integers $u_1, u_2,...,u_r$ such that 
$1=u_1a_1+ u_2a_2+\cdots+u_ra_r$. Then $$h=h^{u_1a_1+u_2a_2+\cdots+u_ra_r}=(h^{a_1})^{u_1}(h^{a_2})^{u_2}\cdots(h^{a_r})^{u_r},$$ 
and $o(h^{a_i})=p_i$ for every $i\in \{1,2,\cdots,r\}$.
This implies that 
$$\langle h^{a_1},h^{a_2},\cdots,h^{a_r}\rangle\geq\langle h\rangle.$$
Clearly, 
$\langle h^{a_1},h^{a_2},\cdots,h^{a_r}\rangle\leq\langle h\rangle$. Thus,
$\langle h^{a_1},h^{a_2},\cdots,h^{a_r}\rangle=\langle h\rangle$,
which means that
$$\langle h\rangle=\bigvee_{1\leq i\leq r}\langle h^{a_i}\rangle.$$
Hence, $h\in \langle t\in H: o(t) \mbox{ is prime}\rangle$, 
which means that $H\leq \langle t\in H: o(t) \mbox{ is prime}\rangle$. 
Clearly, $\langle t\in H: o(t) \mbox{ is prime}\rangle\leq H$. 
Therefore,
$$H=\langle t\in H: o(t) \mbox{ is prime}\rangle,$$ 
which together with Proposition \ref{P02} implies that 
$$H=\bigvee \{A\in A(L(G)): A\leq H\}.$$
Consequently, $L(G)$ is atomistic.
\endproof

\subsection{Cyclic actions}
Let $H$ be a finite group acting linearly on $V$ via a homomorphism
$$\theta:H\longrightarrow \mbox{GL}_{\mathbb{F}_p}(V),$$
and write $$h\cdot v=\theta(h)(v)$$
for all $h\in H$ and $v\in V$. Then $V$ naturally has 
the structure of an $\mathbb{F}_pH$-module. 
The corresponding semidirect product $G=V\rtimes_{\theta} H$ 
is equipped with the multiplication
$$(v,h)(w,k)=(v+h\cdot w,hk)$$
and inverse
$$(v,h)^{-1}=(-h^{-1}\cdot v,h^{-1})$$
for all $v,w\in V$ and $h,k\in H$ \cite{Rotman, Webb}. 
We identify $V$ and $H$ with the subgroups $V\times \{1\}$ 
and $\{0\}\times H$ of $G$, respectively, 
where $1$ denotes the identity element of $H$. 
Accordingly, when no confusion can arise, 
we write $v$ for $(v,1)$ with $v\in V$, and $h$ for $(0,h)$ with $h\in H$.
Moreover, for any $\mathbb{F}_p$-subspace $S\leq V$ and $\rho\in H$, we write 
$$\mbox{Fix}_{S}(\rho)=\{z\in S:\rho\cdot z=z\}.$$
We now consider the special case $H=C_m=\langle b\rangle$. 
Let $\alpha\in \mbox{GL}_{\mathbb{F}_p}(V)$ satisfy $\alpha^m=\mbox{id}_V$. 
We let $C_m$ act on $V$ via $$b\cdot v=\alpha(v),$$
or equivalently, $b^i\cdot v=\alpha^i(v)$. 
We denote the corresponding semidirect product 
by $$G=V\rtimes_{\alpha} C_m.$$
The zero vector of $V$ is denoted by $0_V$, 
or simply by $0$. The same symbol $0$ may 
also denote the zero subspace or zero submodule.

The following lemma follows immediately from Goursat's lemma 
(see, e.g., \cite{Schmidt}). 
\begin{lemma}\label{LeM+1}
Let $H$ and $K$ be finite groups with $\emph{gcd}\{|H|, |K|\}=1$. Then every subgroup $S\leq H\times K$ has a unique expression $S=A\times B$ with $A\leq H$ and $B\leq K$. 
Consequently, $$L(H\times K)\cong L(H)\times L(K).$$   
\end{lemma}

\begin{lemma}\label{L4.1}
Let $m$ be square-free with $p\nmid m$, let $C_m=\langle b\rangle$ and let $\alpha\in \emph{GL}_{\mathbb{F}_p}(V)$ satisfy $\alpha^m=\emph{id}_V$. Put 
$$G=V\rtimes_{\alpha} C_m, d=o(\alpha), e=m/d, $$ and define $$c=b^e, D=\langle c\rangle, N=\langle b^d\rangle, \beta=\alpha^e.$$ 
Let $G_0=\langle V, D\rangle\leq G$. Then $G_0=VD=V\rtimes_\beta D$, the subgroup
$N$ acts trivially on $V$, and the action of $D$ on $V$ given by $c\cdot v=\beta(v)$ is faithful. 
Moreover, \begin{equation}\label{00E1}
G\cong G_0\times N, |D|=d, |N|=e,   
\end{equation}
\begin{equation}\label{00E2}
\mathbb{F}_p[\beta]=\mathbb{F}_p[\alpha]   
\end{equation}
and
\begin{equation}\label{00E3}
L(G)\mbox{ is consistent } \Longleftrightarrow
L(G_0)\mbox{ is consistent}.
\end{equation}
\end{lemma}
\proof 
Since $\alpha^m=\mbox{id}_V$ and $d=o(\alpha)$, 
it follows that $d\mid m$. Because $m$ is square-free 
and $e=m/d$, we know that $\mbox{gcd}\{d,e\}=1$. 
Thus $$|D|=o(b^e)=d \mbox{ and } |N|=o(b^d)=e.$$
Hence $\mbox{gcd}\{|D|,|N|\}=1$, which implies that 
$D\cap N=\{1\}$. Therefore, $$C_m=DN\cong D\times N.$$
Note that $b^d$ acts on $V$ as $\alpha^d=\mbox{id}_V$. 
So, $N$ acts trivially on $V$. Moreover, since $c=b^e$ 
acts on $V$ as $\beta=\alpha^e$ and $\mbox{gcd}\{d,e\}=1$, $o(\beta)=o(\alpha^e)=d=|D|$. 
Therefore, the action of $D=\langle c\rangle$ 
on $V$ given by $c\cdot v=\beta(v)$ is faithful.

As $V\unlhd G$, we have $G_0=VD$. Note that $V\cap D=\{(0,1)\}$. 
Thus we have $$G_0\cong V\rtimes_{\beta}D$$ 
since $D$ acts on $V$ as $\beta$.
Moreover, $N$ centralizes $D$ because $C_m$ is abelian, 
and it also centralizes $V$ since $N$ acts trivially on $V$. 
Thus, $N$ centralizes $G_0$. Furthermore, it is clear that 
$$G=VC_m=VDN=G_0N.$$
We claim that $G_0\cap N=\{(0,1)\}$. Indeed, let $\delta\in G_0\cap N$. 
Since $\delta\in G_0=VD$, there exist $v\in V$ and $y\in D$ such that 
$$\delta=(v,1)(0,y)=(v,y).$$ Moreover, since $\delta\in G_0\cap N\leq C_m$, 
there exists $n\in N$ such that $\delta=(0,n)$. Thus $v=0$ and $y=n$. 
Hence $y\in D\cap N=\{1\}$, and therefore $\delta=(0,1)$.
So, $G_0\cap N=\{(0,1)\}$. Consequently,
$$G\cong G_0\times N$$ since 
$G=G_0N$, $G_0\cap N=\{(0,1)\}$ and $N$ centralizes $G_0$,
which proves \eqref{00E1}.

Clearly, $\mathbb{F}_p[\beta]\subseteq\mathbb{F}_p[\alpha]$. 
We next prove that $\mathbb{F}_p[\alpha]\subseteq\mathbb{F}_p[\beta]$. Indeed, if $d=1$, then $\alpha=\beta=\mbox{id}_V$, 
and $\mathbb{F}_p[\alpha]\subseteq\mathbb{F}_p[\beta]$. Thus \eqref{00E2} is immediate. Now, assume that $d>1$.
Then there exists an integer $1\leq t\leq d-1$ such that $$et\equiv 1 (\mbox{mod }d)$$ because $\mbox{gcd}\{e,d\}=1$. Note that $\beta=\alpha^e$ and $\alpha^d=\mbox{id}_V$. We have 
\begin{equation}\label{+W1}
\beta^t=(\alpha^{e})^t=\alpha^{et}=\alpha.    
\end{equation} 
Thus $\mathbb{F}_p[\alpha]\subseteq\mathbb{F}_p[\beta]$, and hence \eqref{00E2} follows.

Finally, it is clear that $|G_0|=|V||D|=p^nd$. Since $p\nmid m$, $|N|=e$ and $\mbox{gcd}\{d,e\}=1$, we obtain that $\mbox{gcd}\{|G_0|, |N|\}=1$. By Lemma \ref{LeM+1} and \eqref{00E1}, $$L(G)\cong L(G_0)\times L(N).$$ 
Moreover, since $N$ is cyclic, $L(N)$ is modular. Thus $L(N)$ is consistent. Therefore, by Lemma \ref{Lem4.0}, $L(G)$ is consistent if and only if $L(G_0)$ is consistent, which proves \eqref{00E3}.
\endproof

\begin{corollary}\label{Corr4.1}
Under the notation and assumptions of Lemma \ref{L4.1}, the following statements hold.\\
(i1) An $\mathbb{F}_p$-subspace $W\leq V$ is $\alpha$-invariant if and only if it is $\beta$-invariant.\\
(i2) $\mathbb{F}_p[\alpha]$ is a field if and only if $\mathbb{F}_p[\beta]$ is a field.\\
(i3) $m_{\alpha}(x)$ is irreducible over $\mathbb{F}_p$
if and only if $m_{\beta}(x)$ is irreducible over $\mathbb{F}_p$.\\
(i4) Let $C_\alpha=\langle\alpha\rangle$. Then $V$ is
isotypic as an $\mathbb{F}_pC_\alpha$-module if and only if
it is isotypic as an $\mathbb{F}_pD$-module.
\end{corollary}
\proof 
By Lemma \ref{L4.1}, $\mathbb{F}_p[\alpha]=\mathbb{F}_p[\beta]$. Hence, (i2) follows directly. 

We now prove (i1). Let $W\leq V$ be an $\mathbb{F}_p$-subspace. If $W$ is $\alpha$-invariant, then $W$ is $\beta$-invariant since $\beta=\alpha^e$. Conversely, if $W$ is $\beta$-invariant then $W$ is $\alpha$-invariant by \eqref{+W1}. This proves (i1).
 
For (i3), recall that $$\mathbb{F}_p[T]\cong \mathbb{F}_p[x]/(m_T(x))$$
for any $T\in \mbox{End}_{\mathbb{F}_{p}}(V)$. 
Thus, $\mathbb{F}_p[T]$ is a field
if and only if $m_T(x)$ is irreducible over $\mathbb{F}_p$.
Applying this fact to $T=\alpha$ and $T=\beta$, and using (i2), we conclude that (i3) holds.

Finally, we prove (i4). By Lemma \ref{L4.1},
$D=\langle c\rangle$, $c\cdot v=\beta(v)$ and the action of $D$ on $V$ is faithful. 
Moreover, $$\langle \beta\rangle=\langle \alpha^e\rangle=\langle \alpha\rangle=C_\alpha$$ since $\mbox{gcd}\{e,d\}=1$ and $o(\alpha)=d$.
Thus the map 
$$\varphi:D\longrightarrow C_\alpha, \quad \varphi(c^k)=\beta^k$$
is a group isomorphism. 
The action of $D$ on $V$ is precisely 
the action obtained from the natural
action of $C_\alpha$ on $V$ via $\varphi$. Hence the
$\mathbb{F}_pD$-module structure on $V$ and the
$\mathbb{F}_pC_\alpha$-module structure on $V$ correspond under this
isomorphism. Consequently, (i4) follows.
\endproof 

We use Maschke's theorem and Schur's lemma in the following form \cite{Webb, Lux}.
\begin{theorem}[\cite{Webb, Lux}]\label{cq1}
Let $F$ be a field and let $G$ be a finite group
with $\emph{char} F\nmid |G|$. Then every finite-dimensional $FG$-module is semisimple.
\end{theorem}

\begin{lemma}[\cite{Webb,Lux}]\label{LLLL}
Let $F$ be a field and let $G$ be a finite group. 
If $M$ and $N$ are simple $FG$-modules, then every nonzero 
$FG$-homomorphism $f:M\longrightarrow N$ is an isomorphism. 
In particular, if $M\ncong N$, then $Hom_{FG}(M,N)=0$.
\end{lemma}

\begin{proposition}\label{PRO4.1}
Let $G=V\rtimes_{\theta} H$, where $p\nmid |H|$ and $\emph{exp}(H)$ is square-free. Then $L(G)$ is atomistic. Moreover, if $L(G)$ is consistent, then $V$ is isotypic as an $\mathbb{F}_pH$-module.
\end{proposition}
\proof
Since $G/V\cong H$, it follows that $o(gV)\mid \mbox{exp}(H)$ for all $g\in G$. Hence $$(gV)^{\mbox{exp}(H)}=g^{\mbox{exp}(H)}V=V,$$ 
which implies that $g^{\mbox{exp}(H)}\in V$. Note that every element of $V$ has order dividing $p$. Therefore,
$$(g^{\mbox{exp}(H)})^p=g^{\mbox{exp}(H)p}=1.$$
Hence, every $g\in G$ has order dividing $\mbox{exp}(H)p$. Thus,
\begin{equation*}
\mbox{exp}(G)\mid \mbox{exp}(H)p,
\end{equation*}
which, together with $p\nmid |H|$ and the fact that $\mbox{exp}(H)$ is square-free, implies that $\mbox{exp}(G)$ is square-free. Consequently, $L(G)$ is atomistic by Proposition \ref{P03}.

Assume that $L(G)$ is consistent, but that $V$ is not isotypic as an $\mathbb{F}_p H$-module. By Theorem \ref{cq1}, $V$ is semisimple as an $\mathbb{F}_pH$-module since $p\nmid |H|$. 
Then there exist two nonisomorphic simple submodules $U_1, U_2\leq V$. Choose $0\neq u\in U_1$ and $0\neq w\in U_2$, 
and set 
$$I=\langle (u,1)\rangle, K=\langle (w,1)\rangle \mbox{ and } J=\langle (u+w,1)\rangle.$$ 
Since $U_1$ and $U_2$ are simple and nonisomorphic, we have
$U_1\cap U_2=0$. Thus $u+w\neq 0$. 
So, $|I|=|J|=|K|=p$, and hence 
\begin{equation}\label{+EQation}
 I, J, K\in A(L(G))   
\end{equation} 
by Proposition \ref{P02}. 

Since every finite-dimensional $\mathbb{F}_pH$-module is isotypic when $H$ is trivial, 
our assumption that $V$ is not isotypic implies that $H\neq\{1\}$. As $\mbox{exp}(H)$ is square-free, $H$ is generated by its elements of prime order. 
Choose elements $r_1,r_2\cdots r_t\in H$ of prime order such that $H=\langle r_1,r_2\cdots r_t\rangle$, and put 
$$\mathcal{X}=\{\langle(0,r_i)\rangle: 1\leq i\leq t\}.$$
Then $\mathcal{X}\subseteq A(L(G))$
and \begin{equation}\label{G1}
\langle\mathcal{X}\rangle= H.    
\end{equation}
Set $k=u+w$ and $M=\mathbb{F}_{p}H(u+w)=\mathbb{F}_{p}Hk$. 
Then $M$ is an $\mathbb{F}_{p}H$-submodule of $V$, 
so it is invariant under $H$. Thus $M\rtimes_\theta H$ 
is a subgroup of $G$. Note that $J=\langle (k,1)\rangle \leq M\rtimes_\theta H$ 
and $H\leq M\rtimes_\theta H$. Hence, by \eqref{G1}, 
$$\langle\mathcal{X},J\rangle\leq M\rtimes_\theta H.$$
On the other hand, for every $h\in H$, we have 
$(0, h)(k,1)(0,h)^{-1}=(h\cdot k, 1)\in \langle\mathcal{X},J\rangle \cap V$. 
As $\langle\mathcal{X},J\rangle \cap V$ is a subgroup of $V$, it is 
an $\mathbb{F}_{p}$-subspace of $V$. Thus, 
$$M=\mathbb{F}_{p}Hk=\mbox{span}_{\mathbb{F}_p}\{h\cdot k: h\in H\}\leq \langle\mathcal{X},J\rangle \cap V.$$ 
Hence $M\leq \langle\mathcal{X},J\rangle$, 
which together with $H\leq \langle\mathcal{X},J\rangle$ implies that $M\rtimes_\theta H\leq \langle\mathcal{X},J\rangle$. 
Therefore $M\rtimes_\theta H=\langle\mathcal{X},J\rangle$, and hence 
\begin{equation}\label{G2}
\langle\mathcal{X},J\rangle \cap V=\mathbb{F}_{p}Hk=\mathbb{F}_{p}H(u+w). 
\end{equation}

Now, we prove that
\begin{equation}\label{G3}
M=\mathbb{F}_{p}H(u+w)=U_1\oplus U_2.  
\end{equation}
In fact, we have $M\leq U_1\oplus U_2$ since $U_1$ and $U_2$ are $\mathbb{F}_{p}H$-submodules.
Let $$\pi_i:M\longrightarrow U_i,\qquad i=1,2,$$
be the restrictions to $M$ of the canonical projections of $U_1\oplus U_2$ 
onto $U_i$. As $\pi_1(u+w)=u\neq 0$ and $\pi_2(u+w)=w\neq 0$, 
both $\pi_1$ and $\pi_2$ are nonzero $\mathbb{F}_pH$-homomorphisms.
Moreover, both projections are surjective since $U_1$ and $U_2$ are simple. 
Suppose that $M\neq U_1\oplus U_2$. Then $$\ker(\pi_1)=M\cap U_2.$$
We claim that $M\cap U_2$ is a proper submodule of $U_2$. 
Otherwise, $U_2\le M$. Since $\pi_1$ is surjective, 
this would imply $U_1\le M$, and hence $M=U_1\oplus U_2$, a contradiction.
Therefore, since $U_2$ is simple, $\ker(\pi_1)=0$.
Thus $\pi_1:M\to U_1$ is an isomorphism. It follows that
$$\pi_2\pi_1^{-1}:U_1\longrightarrow U_2$$
is a nonzero $\mathbb{F}_pH$-homomorphism, which contradicts 
Lemma \ref{LLLL} because $U_1\ncong U_2$. Consequently, $$M=U_1\oplus U_2=\mathbb{F}_{p}H(u+w).$$

From \eqref{G1}, \eqref{G2} and \eqref{G3}, we get
$\langle\mathcal{X},J\rangle=(U_1\oplus U_2)\rtimes_\theta H$.
Moreover, the same argument as above gives
$\langle X,I,K\rangle=(U_1\oplus U_2)\rtimes_\theta H$.
Hence
\begin{equation}\label{EQ4.33}
\langle\mathcal{X},J\rangle=\langle\mathcal{X},I,K\rangle.\end{equation}
Similar to the proof of \eqref{G2}, we get 
$$\langle\mathcal{X},I\rangle\cap V=\mathbb{F}_{p}Hu=U_1\mbox{ and }\langle\mathcal{X},K\rangle\cap V=\mathbb{F}_{p}Hw=U_2.$$  
By \eqref{EQ4.33}, $J\leq \langle \mathcal{X},I,K\rangle$. Now, we claim that $J\nleq \langle \mathcal{X},I\rangle$. 
Indeed, if $J\leq \langle \mathcal{X},I\rangle$, then $u+w\in \langle \mathcal{X},I\rangle\cap V=\mathbb{F}_{p}Hu=U_1$ since $u+w\in J\cap V$.
It follows that $w=u+w-u\in U_1$, which means that $w\in U_1\cap U_2=0$, contrary to the fact that $w\neq 0$. Similarly, we also have that $J\nleq \langle \mathcal{X},K\rangle$. Thus \begin{equation}\label{+++1}
J\leq \langle \mathcal{X},I,K\rangle,J\nleq \langle \mathcal{X},I\rangle \mbox{ and }J\nleq \langle \mathcal{X},K\rangle, \mbox{ but } \langle \mathcal{X},J\rangle=\langle \mathcal{X},I,K\rangle.    
\end{equation} 
Because $L(G)$ is atomistic, its canonical closure operator on $A(L(G))$ is
$$\mbox{cl}(\mathcal{Y})=\{A\in A(L(G)):A\leq\langle\mathcal{Y}\rangle\}.$$
Therefore, by \eqref{+EQation} and \eqref{+++1},
$$J\in\mbox{cl}(\mathcal{X}\cup\{I,K\})
\setminus [\mbox{cl}(\mathcal{X}\cup\{I\})
\cup \mbox{cl}(\mathcal{X}\cup\{K\})]$$
but
$$\mbox{cl}(\mathcal{X}\cup\{J\})=\mbox{cl}(\mathcal{X}\cup\{I,K\}).$$
Thus the canonical closure operator fails CEP. By Theorem \ref{T04}, $L(G)$ is inconsistent, a contradiction.  
\endproof

\begin{proposition}\label{PRO4.2}
Let $G=V\rtimes_{\theta} H$. 
Suppose that there exist two commuting elements $\rho, \sigma \in H$ of prime order such that $\rho\notin \langle \sigma\rangle$, and a nonzero $\mathbb{F}_p\langle\rho,\sigma\rangle$-submodule
$S\leq V$ such that $$\emph{Fix}_S(\rho)=\emph{Fix}_S(\sigma)=0.$$ Then $L(G)$ is inconsistent.
\end{proposition}
\proof
Set $\rho_S=\theta(\rho)|_S$ and $\sigma_S=\theta(\sigma)|_S$. 
Since $\mbox{Fix}_S(\rho)=0$, the map $\mbox{id}_S-\rho_S$ is injective, and hence 
invertible since $S$ is finite-dimensional. Similarly, the map 
$\mbox{id}_S-\sigma_S$ is also invertible. Choose $0\neq y\in S$ and put 
$$z=(\mbox{id}_S-\rho_S)(y)=y-\rho\cdot y.$$ 
Clearly, $0\neq z\in S$. Moreover,
\begin{equation*}\label{EQ4.3}
(y,1)(0,\rho)(-y,1) =(y+\rho\cdot (-y),\rho) =(y-\rho\cdot y,\rho) =(z,\rho).
\end{equation*}
Hence $(z,\rho)$ is conjugate to $(0,\rho)$, which implies that $o((z,\rho))=o((0,\rho))=o(\rho)$. Define  $J=\langle (z,\rho)\rangle$ and $X=\langle (0,\sigma)\rangle$. Then, by Proposition \ref{P02}, $J,X\in A(L(G))$ since $\rho$ and $\sigma$ have prime order. 

Let $T=\langle X,J\rangle=\langle (0,\sigma), (z,\rho)\rangle$. 
Using the convention $[a, b]=aba^{-1}b^{-1}$, we have 
$$[(z,\rho), (0,\sigma)]=(z,\rho)(0,\sigma)(z,\rho)^{-1}(0,\sigma)^{-1}.$$
Note that $(z,\rho)^{-1}=(-\rho^{-1}\cdot z,\rho^{-1})$. Therefore 
\begin{equation*}
\begin{aligned}
	 [(z,\rho), (0,\sigma)] & =(z,\rho)(0,\sigma)(-\rho^{-1}\cdot z,\rho^{-1}) (0,\sigma)^{-1}\\
& = (z,\rho\sigma)(-\rho^{-1}\cdot z,\rho^{-1}) (0,\sigma^{-1})\\
& = (z+\rho\sigma\cdot (-\rho^{-1}\cdot z),\rho\sigma\rho^{-1})(0,\sigma^{-1})\\
& = (z-\sigma\cdot z,\sigma)(0,\sigma^{-1}) \quad\mbox{ since } \rho \mbox{ and }\sigma \mbox{ commute}\\
& = (z-\sigma\cdot z,1)\\
& = ((\mbox{id}_S-\sigma_S)(z),1).
\end{aligned}    
\end{equation*}
Hence 
\begin{equation}\label{EQ4.4}
[(z,\rho), (0,\sigma)]=((\mbox{id}_S-\sigma_S)(z),1)\in T\cap V.   \end{equation}
Since $(0,\sigma)\in T$, it follows that $(0,\sigma)(w,1)(0,\sigma)^{-1}\in T$ for every $(w,1)\in T\cap V$.
Note that $(0,\sigma)(w,1)(0,\sigma)^{-1}=(\sigma\cdot w,1)\in V$. Thus $(\sigma\cdot w, 1)\in T\cap V$. Hence, $T\cap V$ is an $\mathbb{F}_p$-subspace of $V$, and it is
invariant under $\theta(\sigma)$. Therefore, 
\begin{equation}\label{EQ4.5}
T\cap V  \mbox{ is stable under every element of } \mathbb{F}_p[\theta(\sigma)].
\end{equation} 

Because $\mbox{id}_S-\sigma_S$ is invertible, 
$$\mbox{gcd}\{1-x,m_{\sigma_S}(x)\}=1.$$
Hence, there exist polynomials $a(x), b(x)\in\mathbb{F}_p[x]$ such that
$$a(x)(1-x)+b(x)m_{\sigma_S}(x)=1.$$
Evaluating at $\sigma_S$ gives $a(\sigma_S)(\mbox{id}_S-\sigma_S)=\mbox{id}_S$. Thus $(\mbox{id}_S-\sigma_S)^{-1}=a(\sigma_S)$.
Note that $(\mbox{id}_S-\sigma_S)(z)\in S$. Hence
$$z=a(\sigma_S)(\mbox{id}_S-\sigma_S)(z)=a(\theta(\sigma))(\mbox{id}_S-\sigma_S)(z).$$
From \eqref{EQ4.4} and \eqref{EQ4.5}, we obtain that $$(z,1)\in T\cap V.$$ Moreover, since $(z,\rho)=(z,1)(0, \rho)\in T$, it follows that $$(0, \rho)=(-z,1)(z,\rho)\in T.$$
Similarly, $T\cap V$ is invariant under $\theta(\rho)$ and $(0, \sigma)\in T$.

Let $$K=\langle\rho,\sigma\rangle \mbox { and }W=\mbox{span}_{\mathbb F_p}\{\tau\cdot z:\tau\in K\}.$$
Since $(z,1)\in T\cap V$ and $T\cap V$ is invariant under both $\theta(\rho)$ and $\theta(\sigma)$, we have $W\times\{1\}\leq T$. Furthermore, since $(0,\rho)\in T$ and $(0,\sigma)\in T$, it follows that $\{(0,\tau): \tau\in K\}\leq T$.
Hence $W\rtimes_{\theta} K\leq T$.
The reverse inclusion holds because $T$ is generated by
$(z,\rho)$ and $(0,\sigma)$, and both elements lie in
$W\rtimes_{\theta} K$. Therefore,   
\begin{equation}\label{EQ4.6}
T=W\rtimes_{\theta} K.  
\end{equation}

Next, we define $$M_1=\{(0,\tau): \tau\in K\}\mbox{ and }M_2=\{(w,\tau): w\in W, \tau\in \langle\sigma\rangle\}.$$ 
Since $W$ is invariant under $\theta(\sigma)$, it is invariant under
the action of $\langle\sigma\rangle$. Hence, $M_2$ is a subgroup of $G$. Clearly, $X\leq M_1$ and $X\leq M_2$. Moreover, we know that $M_1<T$ since $(z,1)\in T\setminus M_1$. Also, $M_2<T$ because $(0,\rho)\in T$, but $(0,\rho)\notin M_2$ as $\rho\notin \langle\sigma\rangle$. Finally, since $\{0\}\times K=M_1$ and $W\times\{1\}\leq M_2$, we have $T\geq \langle M_1, M_2\rangle\geq W\rtimes_{\theta} K$. Therefore, by \eqref{EQ4.6},
$$\langle M_1, M_2\rangle=W\rtimes_{\theta} K=T.$$
Consequently, $T=M_1\vee M_2$, where $X\leq M_1, M_2<T$.
Hence, $T=X\vee J$ is not join-irreducible in the interval $[X, G]$. Note that $J\in A(L(G))\subseteq J^{+}(L(G))$. Thus, the join-irreducible element $J$ is inconsistent. Consequently, $L(G)$ is inconsistent.
\endproof
\begin{proposition}[\cite{Webb}]\label{qp03}
Let $F$ be a field and $C_n=\langle c\rangle$. Then $FC_n\cong F[x]/(x^n-1)$ as $F$-algebras, where the residue class of $x$ corresponds to $c$.   
\end{proposition}
\quad\\
\textbf{Proof of Theorem \ref{T001}:}
Retain the notation of Lemma \ref{L4.1}, set $$c=b^e, D=\langle c\rangle, N=\langle b^d\rangle \mbox{ and }\beta=\alpha^e,$$ 
where $e=m/d$, and let $G_0=V\rtimes_\beta D$. The proof is made in two steps as follows. \\
\textbf{Step 1.} Equivalence of (b), (c) and (d).

It suffices to consider the case in which $d$ is prime. 
By Theorem \ref{cq1}, $V$ is semisimple as an $\mathbb{F}_pD$-module since $p\nmid d$. By Proposition \ref{qp03}, 
\begin{equation*}
\mathbb{F}_pD\cong \mathbb{F}_p[x]/(x^d -1).    
\end{equation*}
Moreover, $x^d-1$ is squarefree over $\mathbb{F}_p$ since $(x^d-1)'=dx^{d-1}$ and $p\nmid d$.
Therefore, each simple $\mathbb{F}_pD$-module is precisely the module $\mathbb{F}_p[x]/(f(x))$, where $f(x)$ is 
the monic irreducible factor of $x^d-1$.
Hence there exist monic irreducible factors $f_1(x),\cdots,f_t(x)$ of $x^d-1$, not necessarily distinct, such that
$$V\cong\bigoplus_{i=1}^{t} \mathbb{F}_p[x]/(f_i(x))$$
as an $\mathbb{F}_pD$-module, where the action of the generator $c$ on 
each summand is given by multiplication by the residue class of $x$.
Therefore, $m_\beta(x)=\mbox{lcm}\{f_1(x),\cdots,f_t(x)\}$.
Since every $f_i(x)$ is irreducible, $m_\beta(x)$ is irreducible if and only if all the $f_i(x)$ coincide. 
Consequently, 
\begin{equation*}
V\mbox{ is isotypic}
\Longleftrightarrow m_\beta(x)\mbox{ is irreducible over } \mathbb{F}_p.
\end{equation*}
Note that $$\mathbb{F}_p[\beta]\cong \mathbb{F}_p[x]/(m_\beta(x)).$$
Hence $$\mathbb{F}_p[\beta]\mbox{ is a field }
\Longleftrightarrow
m_\beta(x)\mbox{ is irreducible over }\mathbb{F}_p.$$
By Corollary \ref{Corr4.1}, the corresponding statements for $\alpha$ and
$\beta$ are equivalent. Therefore, $$(b)\Longleftrightarrow(c)\Longleftrightarrow(d).$$\\
\textbf{Step 2.} Equivalence of (a) and (b).

By Lemma \ref{L4.1}, we have $|D|=o(\beta)=d, G\cong G_0\times N$, the action of $D$ on $V$ is faithful, and 
$$L(G) \mbox{ is consistent } \Longleftrightarrow  L(G_0)\mbox{ is consistent }.$$
 Therefore, by Corollary \ref{Corr4.1}, to complete the proof, it is enough to prove that 
$L(G_0)$ is consistent if and only if 
either $d=1$, or $d$ is prime and $V$  is isotypic as 
an $\mathbb{F}_pD$-module. 

Note that since $d\mid m$ and $m$ is squarefree, $\mbox{exp}(D)=d$ is squarefree. Moreover, $p\nmid d=|D|$ since $p\nmid m$. Thus by 
Proposition \ref{PRO4.1}, we know that  
\begin{equation}\label{J+1}
L(G_0) \mbox{ is atomistic}
\end{equation} 
and
\begin{equation}\label{J+21}
 L(G_0) \mbox{ is consistent} \Longrightarrow V \mbox{ is isotypic as an } \mathbb{F}_pD\mbox{-module.}
\end{equation}
\quad\\
\textbf{Necessity.} Suppose that $L(G_0)$ is consistent. 
Assume that $d>1$. 
By \eqref{J+21}, it remains to show that $d$ is prime. 
Suppose for contradiction that $d$ is composite. 
Since $d$ is squarefree, there exist two distinct primes $r, s$ 
such that both $r\mid d$ and $s\mid d$. Put $$\rho=c^{d/r}\mbox{ and }\sigma=c^{d/s}.$$ 
Then $o(\rho)=r$, $o(\sigma)=s$ and $\rho\notin \langle \sigma\rangle$. 
Moreover, $\rho\sigma=\sigma\rho$ since $D$ is cyclic. 
On the other hand, since $V$ is isotypic, there exists a simple $\mathbb{F}_p D$-module $S\leq V$ such that $V\cong S^{\oplus l}$ for an integer $l\geq 1$.
We claim that the action of $D$ on $S$ is faithful. Otherwise, there exists a nonidentity 
element of $D$ acted trivially on $S$. Then the nonidentity 
element of $D$ would act trivially 
on every copy of $S$, and therefore on $V$, contrary to the fact 
that the action of $D$ on $V$ is faithful. 

Since $D$ is abelian, $\mbox{Fix}_{S}(\rho)$ is 
an $\mathbb{F}_pD$-submodule of $S$. Thus 
$$\mbox{Fix}_{S}(\rho)=0 \mbox{ or }\mbox{Fix}_{S}(\rho)=S$$ 
because $S$ is simple. We claim that $\mbox{Fix}_{S}(\rho)\neq S$. 
Otherwise, $\rho$ acts trivially on $S$. 
However, $\rho\neq 1$, which contradicts the 
faithfulness of the action of $D$ on $S$. 
Hence $\mbox{Fix}_{S}(\rho)=0$. 
Similarly, $\mbox{Fix}_{S}(\sigma)=0$. By Proposition \ref{PRO4.2}, $L(G_0)$ is inconsistent, a contradiction. Consequently, $d$ is prime. 

\quad\\
\textbf{Sufficiency.} If $d=1$ then $G_0\cong V$. Thus $G_0$ is an abelian group, which implies that $L(G_0)$ is modular. Hence, $L(G_0)$ is consistent. 

Now, suppose that $d$ is prime and $V$ is isotypic as an $\mathbb F_pD$-module. Denote $F=\mathbb{F}_p[\beta]$. By Step 1, $F$ is a field.
Since $D$ acts faithfully on $V$ and $d>1$, we have $\beta\neq\mbox{id}_V$. Thus $\mbox{id}_V-\beta$ is a nonzero element of $F$ and is therefore invertible. 
Since $\beta^d=\operatorname{id}_V$,  
$$(\mbox{id}_V-\beta)(\mbox{id}_V+\beta+\cdots+\beta^{d-1})=\mbox{id}_V-\beta^d=0.$$ It follows that
\begin{equation}\label{EQ4.7}
\mbox{id}_V+\beta+\cdots+\beta^{d-1}=0. 
\end{equation}
As $F$ is a field, $V$ may naturally be regarded as a vector space over $F$. Now, we claim that 
\begin{equation}\label{EQ4.00}
W\leq V \mbox{ is } \beta\mbox{-invariant} \Longleftrightarrow W\mbox{ is an }F\mbox{-subspace of } V.
\end{equation} Indeed, if $W$ is $\beta$-invariant, then it is stable under every polynomial in $\beta$, and hence under every element of $F$. 
The converse is immediate.

We proceed to describe the atoms of $L(G_0)$, determine the subgroups generated by these atoms, and finally check CEP.

\textbf{(i) The atoms of $L(G_0)$.} 

For $0\neq v\in V$, set $P_v=\langle (v,1) \rangle$. Then $|P_v|=p$ and $P_v=P_{kv}$ for all $k\in \mathbb{F}_p^\times$. 
By Proposition \ref{P02}, $P_v\in A(L(G_0))$. 
For $x\in V$, set $Q_{x}=\langle (x,c) \rangle$. We now prove $Q_{x}\in A(L(G_0))$. By Proposition \ref{P02}, we only need prove $|Q_x|$ is prime. 
By \eqref{EQ4.7} and $c^d=1$, we have 
\begin{equation}\label{eqqq8}(x,c)^d= (x+\beta(x)+\cdots +\beta^{d-1}(x),c^d)=((\mbox{id}_V+\beta+\cdots+\beta^{d-1})(x),1)=(0,1).\end{equation}
Thus $o((x,c))\mid d$.
Let $\pi: G_0\longrightarrow G_0/V$ be the natural quotient map. Then $$\pi((x,c))=(x,c)V=(x,1)(0,c)V=(0,c)V.$$ 
Moreover, the natural isomorphism
$$G_0/V\longrightarrow D,\qquad (x,c^i)V\longmapsto c^i,$$
maps $(0,c)V$ to $c$. Hence $o((0,c)V)=o(c)=d$, which together with
$\pi((x,c))=(0,c)V$ implies that $d=o(\pi((x,c)))\mid o(x,c)$. It follows that $o((x,c))=d$ is prime. 
Let $\mathcal{P}=\{P_v: 0\neq v\in V\}$ and $\mathcal{Q}=\{Q_x: x\in V\}$.
Consequently, $$\mathcal{P}\cup \mathcal{Q}\subseteq A(L(G_0)).$$

We next claim that $A(L(G_0))\subseteq\mathcal{P}\cup \mathcal{Q}$.
Indeed, let $A$ be any atom of $L(G_0)$. Then by Proposition \ref{P02}, $A$ has prime order. Since $|G_0|=p^n d$, $|A|$ is either $p$ or $d$. 
If $|A|=p$, then $A\leq V$. Therefore, $A=P_v$ for some $0\neq v\in V$. If $|A|=d$, then we prove $A=Q_x$ for some $x\in V$. 
Observe that $A\cap V=\{(0,1)\}$ since $\mbox{gcd}(p,d)=1$ and $V$ is a $p$-group. 
Note that $\mbox{ker}(\pi|_A)=A\cap \mbox{ker}\pi = A\cap V$.
Thus, the map $$\pi|_{A}: A\longrightarrow G_0/V$$ is injective. 
By $|A|=d=|G_0/V|$, the map $\pi|_{A}$ is an isomorphism. 
Hence, there exists a unique element $y\in A$ such 
that $\pi(y)=(0,c)V$. Writing $y=(x,c)$ for a unique $x\in V$, we obtain
$A=\langle y \rangle=\langle (x,c)\rangle=Q_x$.
Moreover, $x$ is uniquely determined by $Q_x$. 
Therefore, $A(L(G_0))\subseteq\mathcal{P}\cup \mathcal{Q}$, and hence
\begin{equation}\label{EQ4.8}
A(L(G_0))=\mathcal{P}\cup \mathcal{Q}.
\end{equation}

\textbf{(ii) Subgroups generated by atoms.}

Let $\mathcal{S}\subseteq A(L(G_0))$. If $\mathcal{S}\subseteq \mathcal{P}$, then
\begin{equation}\label{eqqq10} 
\langle \mathcal{S} \rangle=\mbox{span}_{\mathbb{F}_p}\{v\in V:P_v\in \mathcal{S}\}\times \{1\}.
\end{equation}

If $\mathcal{S}\nsubseteq \mathcal{P}$, then $\mathcal{S}\cap \mathcal{Q}\neq \emptyset$ by \eqref{EQ4.8}. Choose $Q_{x_0}\in \mathcal{S}\cap \mathcal{Q}$, and define
\begin{equation*}\label{eqqq11} 
W_{\mathcal{S}}=\mbox{span}_{F}(\{v\in V:P_v\in \mathcal{S}\}\cup \{x-x_0:Q_x\in \mathcal{S}\}).
\end{equation*}
We first verify that $W_{\mathcal{S}}$ is well defined.
If $P_v=P_{v'}$, then $(v',1)\in P_v=\langle(v,1)\rangle$. 
Thus $v'=kv$ for some $k\in\mathbb{F}_p^\times$.
Hence $Fv'=Fv$ since $\mathbb{F}_p\subseteq F$.
Moreover, the parameter $x$ in $Q_x$ is unique by (i). It remains to prove 
that the definition of $W_{\mathcal{S}}$ is independent of the choice of $Q_{x_0}\in\mathcal{S}$. 
Let $Q_{x_1}\in\mathcal{S}$. Then $x-x_1=(x-x_0)-(x_1-x_0)$
for every $Q_{x}\in\mathcal{S}$.
Hence, the $F$-subspace defined using $x_1$ is contained 
in that defined using $x_0$. By symmetry, 
the reverse inclusion also holds.
Consequently, $W_{\mathcal{S}}$ is well defined.

From (\ref{EQ4.00}), $W_{\mathcal{S}}$ is $\beta$-invariant 
since $W_{\mathcal{S}}$ is an $F$-subspace, which together 
with $\beta\in F^\times$ implies 
that $\beta(W_{\mathcal{S}})=W_{\mathcal{S}}$. Thus
$$(x_0,c)(W_{\mathcal{S}}\times\{1\})(x_0,c)^{-1}= W_{\mathcal{S}}\times\{1\}.$$
Hence, $Q_{x_0}$ normalizes $W_{\mathcal{S}}\times\{1\}$.
Moreover, $(W_{\mathcal{S}}\times\{1\})\cap Q_{x_0}
=\{(0,1)\}$ since $Q_{x_0}\cap V=\{(0,1)\}$. 
Therefore, $W_{\mathcal{S}}\rtimes Q_{x_0}$ is a subgroup of $G_0$.
 
We claim that 
\begin{equation}\label{EQ4.9}
\langle \mathcal{S} \rangle=W_{\mathcal{S}}\rtimes Q_{x_0}.    
\end{equation}
Indeed, if $P_v\in \mathcal{S}$, then $(v,1)\in P_v\leq \langle \mathcal{S} \rangle$. 
For any integer $i$, conjugating $(v,1)$ by $(x_0, c)^i$, we obtain that 
$(\beta^i(v),1)\in \langle \mathcal{S} \rangle$. Hence $Fv\times \{1\}\leq \langle \mathcal{S} \rangle\cap V$ since $\langle\mathcal{S}\rangle\cap V$ 
is an $\mathbb{F}_p$-subspace and $F=\mathbb{F}_p[\beta]$. Similarly, if $Q_x\in \mathcal{S}$, then 
$(x,c)(x_0,c)^{-1}=(x-x_0,1)\in  \langle \mathcal{S} \rangle\cap V.$
Again, conjugating $(x-x_0,1)$ by powers of $(x_0,c)$ gives 
$F(x-x_0)\times \{1\}\leq \langle \mathcal{S} \rangle\cap V$. 
Thus $W_{\mathcal{S}}\times \{1\} \leq \langle \mathcal{S} \rangle\cap V.$
Together with $Q_{x_0}\leq\langle\mathcal{S}\rangle$, this gives
$$W_{\mathcal{S}}\rtimes Q_{x_0}\leq\langle\mathcal{S}\rangle.$$
Conversely, every $P_v\in\mathcal{S}$ lies in $W_{\mathcal{S}}\rtimes Q_{x_0}$. 
For every $Q_x\in \mathcal{S}$, $x-x_0\in W_{\mathcal{S}}$. Note that 
$$(x,c)=(x-x_0,1)(x_0,c)\in W_{\mathcal{S}}\rtimes Q_{x_0}.$$
Thus $Q_x\leq W_{\mathcal{S}}\rtimes Q_{x_0}$.
Therefore, every member of $\mathcal{S}$ is contained in $W_{\mathcal{S}}\rtimes Q_{x_0}$, which implies that 
$$\langle\mathcal{S}\rangle\leq W_{\mathcal{S}}\rtimes Q_{x_0}.$$
Consequently, $\langle\mathcal{S}\rangle=W_{\mathcal{S}}\rtimes Q_{x_0}$.

\textbf{(iii) Verification of CEP.} 

In what follows, we focus on proving that the canonical closure operator $\mbox{cl}$ defined on $A(L(G_0))$ satisfies CEP.  

From \eqref{J+1}, we know that $\mbox{cl}(\mathcal{Y}_1)=\mbox{cl}(\mathcal{Y}_2)\Longleftrightarrow \langle \mathcal{Y}_1\rangle=\langle \mathcal{Y}_2\rangle$ 
for any $\mathcal{Y}_1, \mathcal{Y}_2 \subseteq A(L(G_0))$. 
Let $\mathcal{X}\subseteq A(L(G_0))$, and let $I$, $J$ and $K$ be atoms of $L(G_0)$. Suppose that  
\begin{equation}\label{eqqq13}
J\leq \langle \mathcal{X}, I, K\rangle,
\end{equation} 
but  
\begin{equation}\label{eqqq14} 
J\nleq \langle \mathcal{X}, I\rangle \mbox{ and } J\nleq \langle \mathcal{X}, K\rangle.
\end{equation}
It is enough to show that $\langle \mathcal{X}, J\rangle\neq\langle \mathcal{X}, I, K\rangle$.
There are two cases.

Case 1. $\mathcal{X}\cap \mathcal{Q}\neq \emptyset$. 

Choose $Q_{x_0}\in \mathcal{X}\cap \mathcal{Q}$. By \eqref{EQ4.9}, there is an $F$-subspace $W\leq V$ such that 
$\langle \mathcal{X} \rangle=W\rtimes Q_{x_0}$.
For any $F$-subspace $W_1$ with $W\leq W_1\leq V$, we claim that
\begin{equation}\label{+01}
P_v\leq W_1\rtimes Q_{x_0} \Longleftrightarrow v\in W_1     
\end{equation}
and 
\begin{equation}\label{+02}
Q_x\leq W_1\rtimes Q_{x_0} \Longleftrightarrow x-x_0\in W_1 . \end{equation}
Indeed, \eqref{+01} follows from $P_v\leq V$ and $(W_1\rtimes Q_{x_0})\cap V=W_1\times\{1\}$. For \eqref{+02}, if $Q_x\leq W_1\rtimes Q_{x_0}$ 
then $(x,c)\in W_1\rtimes Q_{x_0}$. 
Thus, there exist an element $w\in W_1$ and an integer $1\leq i\leq d$ such 
that $(x,c)=(w,1)(x_0,c)^i$. Note that $(w,1)(x_0,c)^i=(v,c^i)$ for some $v\in V$. 
Hence $i=1$ and $(x,c)=(w,1)(x_0,c)=(w+x_0,c)$, which implies that  
$x-x_0=w\in W_1$. The converse follows immediately from 
$(x,c)=(x-x_0,1)(x_0,c)$. Consequently, \eqref{+02} holds.

Since $W$ is an $F$-subspace of $V$, the quotient $V/W$ is naturally an $F$-vector space. For each atom $A\in A(L(G_0))$, 
define a subspace of $V/W$ by
\begin{equation}\label{eqqq16} 
\Delta(A)=
\begin{cases}
F(v+W), &\mbox{ if }A=P_v,\\
F((x-x_0)+W\bigr), &\mbox{ if }A=Q_x.
\end{cases}
\end{equation}
We first claim that $\Delta(A)$ is well defined. Indeed, if $P_v=P_{v'}$, then
$v'=kv$ for some $k\in\mathbb{F}_p^\times$. Hence $F(v'+W)=F(v+W)$
since $\mathbb{F}_p\subseteq F$. Moreover, the parameter
$x$ in $Q_x$ is unique. By \eqref{EQ4.8}, $\Delta(A)$ is well defined.

By \eqref{+01} and \eqref{+02}, we have $\Delta(A)=0 \Longleftrightarrow  A\leq \langle \mathcal{X} \rangle$. 
More generally, let $A_1, A_2,\cdots,A_k$ be atoms of $L(G_0)$. Using \eqref{EQ4.9}, we obtain that 
$\langle \mathcal{X},A_1, \cdots,A_k\rangle=W'\rtimes Q_{x_0}$
where $W'=W+\sum_{A_i=P_{v_i}}Fv_{i}
+\sum_{A_i=Q_{x_i}}F({x_i-x_0})$.
Consequently,
\begin{equation}\label{eqqq17} 
\frac{W'}{W}=\Delta(A_1)+\cdots \Delta(A_k).
\end{equation}
For every atom $A$ of $L(G_0)$, by \eqref{+01}, \eqref{+02},\eqref{eqqq16} and \eqref{eqqq17}, we obtain that
\begin{equation}\label{eqqq18} 
A\leq \langle \mathcal{X},A_1, \cdots,A_k\rangle \Longleftrightarrow \Delta(A)\leq \sum^{k}_{i=1}\Delta(A_i).
\end{equation}

By \eqref{eqqq13}, \eqref{eqqq14} and \eqref{eqqq18}, $$\Delta(J)\leq \Delta(I)+\Delta(K),$$ while 
\begin{equation}\label{eqqq19}
\Delta(J)\nleq \Delta(I), \mbox{ and }\Delta(J)\nleq \Delta(K).
\end{equation}
Hence $\Delta(J)\neq 0$. We claim that $\Delta(I)\neq 0$. 
Otherwise, $\Delta(J) \leq \Delta(I)+\Delta(K)=\Delta(K)$, 
a contradiction. Similarly, we know that $\Delta(K)\neq 0$.
Thus, $\Delta(J),\Delta(I)$ and $\Delta(K)$ are one-dimensional 
$F$-subspaces of $V/W$. Moreover, we claim $\Delta(I)\neq\Delta(K)$. 
Otherwise, $\Delta(I)+\Delta(K)=\Delta(I)$. This implies $\Delta(J)\leq \Delta(I)$, contrary to \eqref{eqqq19}. 
Thus
$$\mbox{dim}_{F}(\Delta(I)+\Delta(K))=2 \mbox{ whereas } \mbox{dim}_{F}(\Delta(J))=1.$$
Set $\langle \mathcal{X},I,K\rangle\cap V=W_1\times\{1\}$ and 
$\langle \mathcal{X},J\rangle\cap V=W_2\times\{1\}$. 
By (\ref{eqqq17}), we obtain that 
$$\frac{W_1}{W}=\Delta(I)+\Delta(K)\neq \Delta(J)=\frac{W_2}{W},$$ which implies that $\langle \mathcal{X},I,K\rangle\neq \langle \mathcal{X},J\rangle$.

Case 2. $\mathcal{X}\cap \mathcal{Q}=\emptyset$. 

Observe that every atom in $\mathcal{X}$ belongs to $\mathcal{P}$ by \eqref{EQ4.8}. 
Put $U=\mbox{span}_{\mathbb{F}_p}\{v\in V: P_v\in \mathcal{X}\}$. 
Then $\langle \mathcal{X}\rangle=U\times\{1\}$ by \eqref{eqqq10}. We consider the following two subcases.

Subcase 1. Suppose that $I, K\in \mathcal{P}$. Then $I=P_i$ and $K=P_k$ for some $0\neq i,k\in V$. 
Thus $$\langle \mathcal{X},I,K\rangle=(U+\mathbb{F}_{p}i+\mathbb{F}_{p}k)\times \{1\}.$$
By \eqref{eqqq13}, we have $J=P_j\in \mathcal{P}$ for some $0\neq j\in V$.
Hence, we know that $\langle \mathcal{X},J\rangle=(U+\mathbb{F}_{p}j)\times \{1\}$. In the quotient 
space $V/U$, define $$\varGamma(P_v)=\mathbb{F}_{p}(v+U)\leq V/U$$ 
for each atom  $P_v\in \mathcal{P}$. 
As in Case 1, \eqref{eqqq13} and \eqref{eqqq14} imply that
$$\varGamma(J)\leq \varGamma(I)+\varGamma(K), 
\mbox{ but }\varGamma(J)\nleq \varGamma(K) \mbox{ and }\varGamma(J)\nleq \varGamma(I).$$
Exactly as in Case 1, we know that $\varGamma(I), \varGamma(K)$ and $\varGamma(J)$ 
are distinct one-dimensional $\mathbb{F}_{p}$-subspaces of $V/U$. Consequently, $$\mbox{dim}_{\mathbb{F}_{p}}(\varGamma(I)+\varGamma(K))=2 \neq 1=\mbox{dim}_{\mathbb{F}_{p}}(\varGamma(J)).$$
Since $$\frac{\langle X,I,K\rangle\cap V}{U\times\{1\}}\cong\Gamma(I)+\Gamma(K) \mbox{ and }
\frac{\langle X,J\rangle\cap V}{U\times\{1\}}\cong\Gamma(J),$$ we have $\langle \mathcal{X},I,K\rangle\neq \langle \mathcal{X},J\rangle$.

Subcase 2. Suppose that $I\in \mathcal{Q}$ or $K\in \mathcal{Q}$. By symmetry, assume that $I=Q_x$ for some $x\in V$. 
Then clearly $\langle \mathcal{X},I,K\rangle\nleq V$. 
If $J\in \mathcal{P}$ then $\langle \mathcal{X},J\rangle\leq V$. Thus $\langle \mathcal{X},I,K\rangle\neq \langle \mathcal{X},J\rangle$. We may now assume that $J\notin \mathcal{P}$. Then
$J\in \mathcal{Q}$ by \eqref{EQ4.8}. Set $J=Q_y$ with $y\in V$. 
Let $W=\mbox{span}_{F}(U)$. By \eqref{EQ4.9}, 
$\langle \mathcal{X},I\rangle=W\rtimes Q_x$.   
Put $\langle \mathcal{X},I, K\rangle\cap V=W_1\times \{1\}$. We prove that 
\begin{equation}\label{HE+}W_1>W.\end{equation} 
In fact, if $K\notin \mathcal{Q}$ then $K=P_v\in \mathcal{P}$ 
for some $0\neq v\in V$. 
We first observe that $K\nleq \langle \mathcal{X},I\rangle.$
Indeed, if $K\leq \langle \mathcal{X},I\rangle$, then
$\langle \mathcal{X},I,K\rangle=\langle \mathcal{X},I\rangle$.
Together with $J\leq\langle \mathcal{X},I,K\rangle$, this would imply
$J\leq\langle\mathcal{X},I\rangle$, contrary to \eqref{eqqq14}.
From $K\nleq \langle \mathcal{X},I\rangle$, we obtain that $v\notin W$.
By \eqref{EQ4.9}, 
\begin{equation}\label{EQ4.11}
W_1=W+Fv>W.   
\end{equation}
If $K\in \mathcal{Q}$ then $K=Q_z$ for some $z\in V$. As above, 
$K\nleq \langle \mathcal{X},I\rangle$, which implies that $z-x\notin W$.
Hence, by \eqref{EQ4.9}, 
\begin{equation}\label{EQ4.12}
W_1=W+F(z-x)>W.   
\end{equation}
Equations \eqref{EQ4.11} and \eqref{EQ4.12} imply \eqref{HE+}. 
Note that $J=Q_y$. Then, by \eqref{EQ4.9}, we have 
$$\langle \mathcal{X},J\rangle\cap V=W\times\{1\},$$ which together with \eqref{HE+} 
implies that $\langle \mathcal{X},I,K\rangle\neq \langle \mathcal{X},J\rangle$.

By Cases 1 and 2, the canonical closure operator $\mbox{cl}$ satisfies CEP. Consequently, by Theorem \ref{T04}, 
$L(G_0)$ is consistent. 

From Steps 1 and 2, the proof is completed. 
\endproof

\begin{example}\label{example001}
\emph{Let $V=\mathbb{F}_{3}e\oplus \mathbb{F}_{3}f$ and $C_2=\langle b\rangle$. Relative to the ordered basis $\{e,f\}$, define $\alpha\in \mbox{GL}_{\mathbb{F}_{3}}(V)$ by
$$\alpha=\begin{bmatrix}
	    1 & 0 \\
	    0 & -1
	\end{bmatrix}.$$
Then $V\cong C^{2}_{3}$ and $o(\alpha)=2$. 
Moreover, $$V=\{ae+cf: a, c\in \mathbb{F}_{3}\}$$ and 
$$\alpha(ae+cf)=ae-cf$$ for all $a,c\in \mathbb{F}_{3}$. 
Let $G=V\rtimes_{\alpha} C_2$. Then $G=\{(ae+cf,b^{i}): a,c\in \mathbb{F}_{3}, i\in \{0,1\}\}$ and $|G|=|V|\cdot |C_2|=9\cdot 2=18$.
Set $U_{1}=\mathbb{F}_{3}e$ and $U_{2}=\mathbb{F}_{3}f$. Then both $U_1$ and $U_2$ 
are simple $\mathbb{F}_{3}C_2$-submodules of $V$. Since $\alpha(e)=e$ and $\alpha(f)=-f$, we have 
$$U_1\cong \mathbb{F}_{3}[x]/(x-1)$$
and 
$$U_2\cong \mathbb{F}_{3}[x]/(x+1)$$
as $\mathbb{F}_3C_2$-modules.
These two simple modules are not isomorphic because $b$ acts as $\mbox{id}_{U_1}$ on $U_1$
and as $-\mbox{id}_{U_2}$ on $U_2$. Hence, $V=U_1\oplus U_2$ is not an isotypic $\mathbb{F}_{3}C_2$-module.
Therefore, by Theorem \ref{T001}, the lattice $L(G)$ is inconsistent.}

\emph{We now describe the subgroup lattice $L(G)$ explicitly. For convenience, set $$\langle (0,b)\rangle=Q_0, \langle (f,b)\rangle=Q_1, \langle (2f,b)\rangle=Q_2,$$ 
$$\langle (e,1)\rangle=P_1, \langle (f,1)\rangle=P_2, \langle (e+f,1)\rangle=P_3, \langle (e-f, 1)\rangle=P_4$$ 
and 
$$\langle (e,1),(0,b)\rangle=C_1, \langle (e,1),(f,b)\rangle=C_2, \langle (e,1), (2f,b)\rangle=C_3, \langle (f,1),(0,b)\rangle=C_4.$$ The complete list of subgroups of $G$ is presented in Table \ref{Tab:01}. It is straightforward to verify that $L(G)$ is depicted in Figure 10.}
\begin{table}[!h]
\centering
\caption{All subgroups of $G$.}
\label{Tab:01}
\begin{tabular}{c|c}
order & subgroups \\
 \hline
 $1$ & $\langle (0,1)\rangle$  \\
 $2$ & $Q_i \quad (0\leq i\leq 2)$  \\
 $3$ & $P_i \quad (1\leq i\leq 4)$ \\
 $6$ & $C_i \quad (1\leq i\leq 4)$\\
 $9$ & $V$ \\
 $18$ & $G$\\
 \end{tabular}
\end{table}

\begin{center}
\begin{minipage}{11pc}
\setlength{\unitlength}{1pt}\begin{picture}(600,180)
\put(60,40){\circle{6}}\put(48,30){\makebox(0,0)[l]{\footnotesize $\langle (0,1)\rangle$}}
\put(100,80){\circle{6}} \put(105,80){\makebox(0,0)[l]{\footnotesize$Q_0$}}
\put(140,80){\circle{6}} \put(145,80){\makebox(0,0)[l]{\footnotesize$Q_1$}}
\put(180,80){\circle{6}} \put(184,80){\makebox(0,0)[l]{\footnotesize$Q_2$}}
\put(20,80){\circle{6}}\put(4,80){\makebox(0,0)[l]{\footnotesize$P_1$}}
\put(60,80){\circle{6}}\put(43,80){\makebox(0,0)[l]{\footnotesize$P_2$}}
\put(-20,80){\circle{6}}\put(-35,82){\makebox(0,0)[l]{\footnotesize$P_3$}}
\put(-60,80){\circle{6}}\put(-77,80){\makebox(0,0)[l]{\footnotesize$P_4$}}
\put(-20,120){\circle{6}}\put(-35,120){\makebox(0,0)[l]{\footnotesize$V$}}
\put(20,120){\circle{6}}\put(5,120){\makebox(0,0)[l]{\footnotesize$C_1$}}
\put(60,120){\circle{6}}\put(47,128){\makebox(0,0)[l]{\footnotesize$C_2$}}
\put(100,120){\circle{6}}\put(105,122){\makebox(0,0)[l]{\footnotesize $C_3$}}
\put(140,120){\circle{6}}\put(145,120){\makebox(0,0)[l]{\footnotesize $C_4$}}
\put(60,160){\circle{6}}\put(55,168){\makebox(0,0)[l]{\footnotesize$G$}}
\put(58,42){\line(-1,1){36}}
\put(58,42){\line(-2,1){75}}
\put(58,42){\line(-3,1){115}}
\put(62,42){\line(2,1){75}}
\put(62,42){\line(3,1){115}}
\put(20,83){\line(0,1){34}}
\put(62,42){\line(1,1){36}}
\put(-20,83){\line(0,1){34}}
\put(-58,82){\line(1,1){36}}
\put(18,82){\line(-1,1){36}}
\put(58,82){\line(-2,1){75}}
\put(22,82){\line(2,1){75}}
\put(98,82){\line(-2,1){75}}
\put(62,82){\line(2,1){75}}
\put(138,82){\line(-2,1){75}}
\put(178,82){\line(-2,1){75}}
\put(178,82){\line(-1,1){36}}
\put(140,83){\line(0,1){34}}
\put(102,82){\line(1,1){36}}
\put(22,82){\line(1,1){36}}
\put(60,123){\line(0,1){34}}
\put(22,122){\line(1,1){36}}
\put(-18,122){\line(2,1){75}}
\put(98,122){\line(-1,1){36}}
\put(138,122){\line(-2,1){75}}
\put(60,43){\line(0,1){34}}
\put(-10,5){\emph{Figure 10. The lattice $L(G)$}}
\end{picture}
\end{minipage}
\end{center}
\emph{Figure 10 shows that $P_3\in J(L(G))$. However, $P_3$ is inconsistent since $P_3\vee Q_0=G$ but $G$ is not a join-irreducible element in the sublattice $[Q_0, G]$.} 
\end{example}

\begin{example}\label{example002}
\emph{Let $V=\mathbb{F}_{7}$ and $C_6=\langle c\rangle$. Define $\alpha\in \mbox{GL}_{\mathbb{F}_{7}}(V)$ by $$\alpha(x)=3x$$ for all $x\in V$.
It is clear that $o(\alpha)=6$. Moreover, $V$ is isotypic as 
an $\mathbb{F}_7C_6$-module and $m_{\alpha}(x)=x-3$ is irreducible over $\mathbb{F}_{7}$. 
Let $G=V\rtimes_{\alpha} C_6$.  By Theorem \ref{T001}, $L(G)$ is inconsistent.}

\emph{We now describe the subgroup lattice $L(G)$ explicitly. For convenience, we denote $$A_i=\langle (i,c)^3\rangle, B_i=\langle (i,c)^2\rangle \mbox{ and } Q_i=\langle (i, c)\rangle \mbox{ for all } 0\leq i\leq 6.$$ 
The complete list of all subgroups of $G$ is presented in Table \ref{Tab:02}. The subgroup lattice $L(G)$ is depicted in Figure 11.}
\begin{table}[!h]
\centering
\caption{All subgroups of $G$.}
\label{Tab:02}
\begin{tabular}{c|c}
order & subgroups \\
 \hline
 $1$ & $\langle (0,1)\rangle$ \\
 $2$ & $A_i \quad (0\leq i\leq 6)$  \\
 $3$ & $B_i \quad (0\leq i\leq 6)$ \\
 $6$ & $Q_i \quad (0\leq i\leq 6)$  \\
 $7$ & $V$  \\
 $14$ & $V\rtimes_{\alpha^3} \langle c^3\rangle$ \\
 $21$ & $V\rtimes_{\alpha^2} \langle c^2\rangle$ \\
 $42$ & $G$ \\
 \end{tabular}
\end{table}

\begin{center}
\begin{minipage}{13pc}
\setlength{\unitlength}{1pt}\begin{picture}(600,180)
\put(60,40){\circle{6}}\put(50,32){\makebox(0,0)[l]{\tiny$\langle (0,1)\rangle$}}
\put(80,80){\circle{6}}\put(83,80){\makebox(0,0)[l]{\tiny$B_0$}}
\put(100,80){\circle{6}} \put(103,80){\makebox(0,0)[l]{\tiny$B_1$}}
\put(120,80){\circle{6}} \put(123,80){\makebox(0,0)[l]{\tiny$B_2$}}
\put(140,80){\circle{6}} \put(143,80){\makebox(0,0)[l]{\tiny$B_3$}}
\put(160,80){\circle{6}} \put(163,80){\makebox(0,0)[l]{\tiny$B_4$}}
\put(180,80){\circle{6}} \put(183,80){\makebox(0,0)[l]{\tiny$B_5$}}
\put(220,80){\circle{6}} \put(223,80){\makebox(0,0)[l]{\tiny$B_6$}}
\put(40,80){\circle{6}} \put(27.5,80){\makebox(0,0)[l]{\tiny$A_6$}}
\put(20,80){\circle{6}}\put(7.5,80){\makebox(0,0)[l]{\tiny$A_5$}}
\put(60,80){\circle{6}}\put(51,79){\makebox(0,0)[l]{\tiny$V$}}
\put(0,80){\circle{6}} \put(-12.5,80){\makebox(0,0)[l]{\tiny$A_4$}}
\put(-20,80){\circle{6}}\put(-32.5,80){\makebox(0,0)[l]{\tiny$A_3$}}
\put(-40,80){\circle{6}} \put(-52.5,80){\makebox(0,0)[l]{\tiny$A_2$}}
\put(-60,80){\circle{6}}\put(-72,80){\makebox(0,0)[l]{\tiny$A_1$}}
\put(-100,80){\circle{6}} \put(-113,80){\makebox(0,0)[l]{\tiny$A_0$}}
\put(-20,120){\circle{6}}\put(-63,123){\makebox(0,0)[l]{\tiny$V\rtimes_{\alpha^3} \langle c^3\rangle$}}
\put(140,120){\circle{6}}\put(144,122){\makebox(0,0)[l]{\tiny$V\rtimes_{\alpha^2} \langle c^2\rangle$}}
\put(20,120){\circle{6}}\put(23,121){\makebox(0,0)[l]{\tiny$Q_1$}}
\put(60,120){\circle{6}}\put(63,121){\makebox(0,0)[l]{\tiny$Q_3$}}
\put(80,120){\circle{6}}\put(83,121){\makebox(0,0)[l]{\tiny$Q_4$}}
\put(100,120){\circle{6}}\put(103,121){\makebox(0,0)[l]{\tiny$Q_5$}}
\put(120,120){\circle{6}}\put(123,121){\makebox(0,0)[l]{\tiny$Q_6$}}
\put(40,120){\circle{6}}\put(43,121){\makebox(0,0)[l]{\tiny$Q_2$}}
\put(0,120){\circle{6}}\put(3,120){\makebox(0,0)[l]{\tiny$Q_0$}}
\put(60,160){\circle{6}}\put(55,168){\makebox(0,0)[l]{\tiny$G$}}
\put(58,42){\line(-1,1){36}}
\put(58,42){\line(-1,2){17.5}}
\put(58,42){\line(-3,2){55}}
\put(57.5,41.5){\line(-2,1){74.5}}
\put(57,41){\line(-3,1){114.5}}
\put(57,41){\line(-5,2){94.5}}
\put(57,41){\line(-4,1){154.5}}
\put(62.5,41.5){\line(2,1){74.5}}
\put(62.5,41.5){\line(3,1){114.5}}
\put(62,42){\line(1,2){17.5}}
\put(62,42){\line(3,2){55}}
\put(63,41){\line(5,2){94.5}}
\put(63,41){\line(4,1){154.5}}
\put(-20,83){\line(0,1){34}}
\put(140,83){\line(0,1){34}}
\put(60,43){\line(0,1){34}}
\put(60,123){\line(0,1){34}}
\put(-2,82){\line(-1,2){17.5}}
\put(38,82){\line(-3,2){55}}
\put(-38,82){\line(1,2){17.5}}
\put(-98.5,81.5){\line(2,1){75.5}}
\put(82,82){\line(3,2){55}}
\put(122,82){\line(1,2){17.5}}
\put(158,82){\line(-1,2){17.5}}
\put(218.5,81.5){\line(-2,1){75.5}}
\put(-97,81){\line(5,2){94.5}}
\put(-57.5,81.5){\line(2,1){74.5}}
\put(-37.5,81.5){\line(2,1){74.5}}
\put(-17.5,81.5){\line(2,1){74.5}}
\put(2.5,81.5){\line(2,1){74.5}}
\put(22.5,81.5){\line(2,1){74.5}}
\put(42.5,81.5){\line(2,1){74.5}}
\put(77.5,81.5){\line(-2,1){74.5}}
\put(97.5,81.5){\line(-2,1){74.5}}
\put(117.5,81.5){\line(-2,1){74.5}}
\put(137.5,81.5){\line(-2,1){74.5}}
\put(157.5,81.5){\line(-2,1){74.5}}
\put(217,81){\line(-5,2){94.5}}
\put(177.5,81.5){\line(-2,1){74.5}}
\put(62,42){\line(1,1){36}}
\put(-58,82){\line(1,1){36}}
\put(18,82){\line(-1,1){36}}
\put(58,82){\line(-2,1){75}}
\put(62,82){\line(2,1){75}}
\put(178,82){\line(-1,1){36}}
\put(102,82){\line(1,1){36}}
\put(22,122){\line(1,1){36}}
\put(-18,122){\line(2,1){75}}
\put(98,122){\line(-1,1){36}}
\put(138,122){\line(-2,1){75}}
\put(2,122){\line(3,2){55}}
\put(42,122){\line(1,2){17.5}}
\put(78,122){\line(-1,2){17.5}}
\put(118,122){\line(-3,2){55}}
\put(-10,5){\emph{Figure 11. The lattice $L(G)$}}
\end{picture}
\end{minipage}
\end{center}
\emph{Figure 11 shows that $A_3\in J(L(G))$ and $A_3\vee B_0=G$. Moreover, we have that $B_0\leq Q_0<G$, $B_0\leq V\rtimes_{\alpha^2}\langle c^2\rangle<G$ and $G=Q_0\vee (V\rtimes_{\alpha^2} \langle c^2\rangle)$. Thus, $G$ is not a join-irreducible element in the sublattice $[B_0, G]$, which implies that $L(G)$ is inconsistent.}
\end{example}

\begin{example}\label{example003}
\emph{Let $V=\mathbb{F}_{7}$ and $C_6=\langle c\rangle$. Define $\alpha\in \mbox{GL}_{\mathbb{F}_{7}}(V)$ by $$\alpha(x)=2x$$ for any $x\in V$.
Then $o(\alpha)=3$. Moreover, $V$ is isotypic as an $\mathbb{F}_7C_6$-module. 
Let $G=V\rtimes_{\alpha} C_6$. Because $\alpha^3=\mbox{id}_V$, the subgroup $N=\langle c^3\rangle$ acts trivially on $V$. 
On the other hand, the subgroup $D=\langle c^2\rangle$ acts faithfully on $V$ via $\beta$, where $\beta=\alpha^2$. 
Clearly, $o(\beta)=3$. Put $H=V\rtimes_{\beta} \langle c^2\rangle$. Then it is clear that $G\cong H\times N$. By Theorem \ref{T001}, $L(G)$ is consistent.} 

\emph{For convenience, we denote $$B_i=\langle (i,c)^2\rangle=\langle (3i,c^2)\rangle \mbox{ and } Q_i=\langle (i, c)\rangle \mbox{ for all } 0\leq i\leq 6.$$ Then the complete list of all subgroups of $G$ is presented in Table \ref{Tab:03}. Thus, $L(G)$ is depicted in Figure 12. 
It can be verified that every element in $J^{+}(L(G))$ is consistent, and hence $L(G)$ is consistent.}
\begin{table}[!h]
\centering
\caption{All subgroups of $G$.}
\label{Tab:03}
\begin{tabular}{c|c}
order & subgroups \\
 \hline
 $1$ & $\langle (0,1)\rangle$ \\
 $2$ & $N$  \\
 $3$ & $B_i \quad (0\leq i\leq 6)$ \\
 $6$ & $Q_i \quad (0\leq i\leq 6)$  \\
 $7$ & $V$  \\
 $14$ & $V\times N$ \\
 $21$ & $H$ \\
 $42$ & $G$ \\
 \end{tabular}
\end{table}

\begin{center}
\begin{minipage}{13pc}
\setlength{\unitlength}{1pt}\begin{picture}(600,180)
\put(60,40){\circle{6}}\put(50,32){\makebox(0,0)[l]{\tiny$\langle (0,1)\rangle$}}
\put(80,80){\circle{6}}\put(83,80){\makebox(0,0)[l]{\tiny$B_4$}}
\put(100,80){\circle{6}} \put(103,80){\makebox(0,0)[l]{\tiny$B_5$}}
\put(120,80){\circle{6}} \put(123,80){\makebox(0,0)[l]{\tiny$B_6$}}
\put(140,80){\circle{6}} \put(143,80){\makebox(0,0)[l]{\tiny$V$}}
\put(40,80){\circle{6}} \put(27.5,80){\makebox(0,0)[l]{\tiny$B_2$}}
\put(20,80){\circle{6}}\put(7.5,80){\makebox(0,0)[l]{\tiny$B_1$}}
\put(60,80){\circle{6}}\put(48,79){\makebox(0,0)[l]{\tiny$B_3$}}
\put(0,80){\circle{6}} \put(-12.5,80){\makebox(0,0)[l]{\tiny$B_0$}}
\put(-20,80){\circle{6}}\put(-31,80){\makebox(0,0)[l]{\tiny$N$}}
\put(-20,120){\circle{6}}\put(-46,121){\makebox(0,0)[l]{\tiny$V\times N$}}
\put(140,120){\circle{6}}\put(144,122){\makebox(0,0)[l]{\tiny$H$}}
\put(20,120){\circle{6}}\put(23,121){\makebox(0,0)[l]{\tiny$Q_1$}}
\put(60,120){\circle{6}}\put(63,121){\makebox(0,0)[l]{\tiny$Q_3$}}
\put(80,120){\circle{6}}\put(83,121){\makebox(0,0)[l]{\tiny$Q_4$}}
\put(100,120){\circle{6}}\put(103,121){\makebox(0,0)[l]{\tiny$Q_5$}}
\put(120,120){\circle{6}}\put(123,121){\makebox(0,0)[l]{\tiny$Q_6$}}
\put(40,120){\circle{6}}\put(43,121){\makebox(0,0)[l]{\tiny$Q_2$}}
\put(0,120){\circle{6}}\put(3,120){\makebox(0,0)[l]{\tiny$Q_0$}}
\put(60,160){\circle{6}}\put(55,168){\makebox(0,0)[l]{\tiny$G$}}
\put(58,42){\line(-1,1){36}}
\put(58,42){\line(-1,2){17.5}}
\put(58,42){\line(-3,2){55}}
\put(57.5,41.5){\line(-2,1){74.5}}
\put(62.5,41.5){\line(2,1){74.5}}
\put(62,42){\line(1,2){17.5}}
\put(62,42){\line(3,2){55}}
\put(0,83){\line(0,1){34}}
\put(20,83){\line(0,1){34}}
\put(40,83){\line(0,1){34}}
\put(60,83){\line(0,1){34}}
\put(80,83){\line(0,1){34}}
\put(100,83){\line(0,1){34}}
\put(120,83){\line(0,1){34}}
\put(-20,83){\line(0,1){34}}
\put(140,83){\line(0,1){34}}
\put(60,43){\line(0,1){34}}
\put(60,123){\line(0,1){34}}
\put(122,82){\line(1,2){17.5}}
\put(-17.5,81.5){\line(2,1){74.5}}
\put(-18,82){\line(1,2){17.5}}
\put(-18,82){\line(1,1){36}}
\put(-17.5,81.5){\line(3,2){54.5}}
\put(-17.5,81.5){\line(5,2){94.5}}
\put(-17.5,81.5){\line(3,1){114.5}}
\put(-17.5,81.5){\line(7,2){134.5}}
\put(62.5,81.5){\line(2,1){74.5}}
\put(102,82){\line(1,1){36}}
\put(82.5,81.5){\line(3,2){54.5}}
\put(42.5,81.5){\line(5,2){94.5}}
\put(22.5,81.5){\line(3,1){114.5}}
\put(2.5,81.5){\line(7,2){134.5}}
\put(137.5,80){\line(-4,1){154.5}}
\put(62,42){\line(1,1){36}}
\put(22,122){\line(1,1){36}}
\put(-18,122){\line(2,1){75}}
\put(98,122){\line(-1,1){36}}
\put(138,122){\line(-2,1){75}}
\put(2,122){\line(3,2){55}}
\put(42,122){\line(1,2){17.5}}
\put(78,122){\line(-1,2){17.5}}
\put(118,122){\line(-3,2){55}}
\put(-10,5){\emph{Figure 12. The lattice $L(G)$}}
\end{picture}
\end{minipage}
\end{center}
\end{example}

\subsection{General faithful actions}
In the following, we shall study arbitrary faithful actions of finite groups.
\begin{proposition}\label{PPOO}
Let $G=V\rtimes_\theta H$ where $H$ acts faithfully on $V$, 
$p\nmid |H|$ and $\emph{exp}(H)$ is square-free. If $L(G)$ is consistent then the following statements hold.\\
(a1) $L(H)$ is consistent.\\
(a2) $V$ is isotypic as an $\mathbb{F}_pH$-module.\\
(a3) If $\rho, \sigma\in H$ are commuting elements of prime order such that $\rho\notin \langle\sigma\rangle$ and $0\neq S\leq V$ is an $\mathbb{F}_p\langle \rho, \sigma\rangle$-submodule, then 
$$\emph{Fix}_{S}(\rho)\neq 0 \mbox{ or } \emph{Fix}_{S}(\sigma)\neq 0.$$
\end{proposition}
\proof
Since $p\nmid |H|$ and $\mbox{exp}(H)$ is square-free, Proposition \ref{P03} implies that $L(H)$ is atomistic. We first prove (a1). Suppose, to the contrary, that $L(H)$ is inconsistent. 
By Theorem \ref{T04}, the canonical closure operator of $L(H)$ fails CEP. Hence, there exist atoms $I, J, K \in A(L(H))$ and $\mathcal{X}\subseteq A(L(H))$ such that
$$J\leq \langle \mathcal{X}, I, K\rangle_{H}, J\nleq \langle \mathcal{X}, I\rangle_{H}, J\nleq \langle \mathcal{X}, K\rangle_{H}$$
and
$$\langle \mathcal{X}, J\rangle_{H}=\langle \mathcal{X}, I, K\rangle_{H}.$$
Note that $H\leq G$. For every family $\mathcal{Y}$ of subgroups of $H$,
$$\langle\mathcal{Y}\rangle_H=
\langle\bigcup_{Y\in\mathcal{Y}}Y\rangle_H=
\langle\bigcup_{Y\in\mathcal{Y}}Y\rangle_G=\langle\mathcal{Y}\rangle_G$$
since the subgroup of $G$ generated by $\mathcal{Y}$ is contained in $H$.
Moreover, we claim that $$A(L(H))\subseteq A(L(G)).$$ 
Indeed, every atom of $L(H)$ has prime 
order, and hence is also an atom of $L(G)$ by Proposition \ref{P02}. 
Therefore, $I,J,K\in A(L(G))$ and $\mathcal{X}\subseteq A(L(G))$, and they satisfy
$$J\leq \langle \mathcal{X}, I, K\rangle_{G}, J\nleq \langle \mathcal{X}, I\rangle_{G}, J\nleq \langle \mathcal{X}, K\rangle_{G},$$
and
$$\langle \mathcal{X}, J\rangle_{G}=\langle \mathcal{X}, I, K\rangle_{G}$$
since $\langle\mathcal{Y}\rangle_{H}= \langle\mathcal{Y}\rangle_G$.
By Proposition \ref{PRO4.1}, $L(G)$ is atomistic. Hence, the canonical closure operator of $L(G)$ also fails CEP. Consequently, by Theorem \ref{T04}, $L(G)$ is inconsistent, yielding a contradiction. Thus, (a1) is proved. 

Statement (a2) follows immediately from Proposition \ref{PRO4.1}. Thus, it remains to prove (a3). Suppose, to the contrary, that there exist
commuting elements $\rho,\sigma\in H$ of prime order with $\rho\notin\langle\sigma\rangle$,
and a nonzero $\mathbb F_p\langle\rho,\sigma\rangle$-submodule $S\leq V$ 
such that both $\mbox{Fix}_{S}(\rho)=0$ and $\mbox{Fix}_{S}(\sigma)=0$. 
Then Proposition \ref{PRO4.2} implies that $L(G)$ is inconsistent, a contradiction. Hence (a3) follows.
\endproof
\begin{example}\label{EX++}
\emph{Let $V=\{(x_1, x_2, x_3, x_4)\in \mathbb{F}^{4}_{7}: x_1+x_2+x_3+x_4=0\}$. Then $V$ is a three-dimensional $\mathbb{F}_7$-vector space, and hence $(V, +)\cong C^3_{7}$. Let the alternating group $A_4$ act naturally on $V$ by permuting the four coordinates. Since $V$ is invariant under this action, it induces a faithful action $$\theta: A_4\longrightarrow \mbox{GL}_{\mathbb{F}_{7}}(V).$$
Moreover, it is straightforward to verify that $V$ is a simple
$\mathbb{F}_{7}A_4$-module. Thus $V$ is isotypic as an $\mathbb{F}_{7}A_4$-module. Let $G=V\rtimes_\theta A_4$. Note that $7\nmid 12=|A_4|$ and $\mbox{exp}(A_4)=6$ is square-free. 
Let $$\rho=(12)(34)\mbox{ and } \sigma=(13)(24)\in A_4.$$ 
Then $\rho$ and $\sigma$ are distinct commuting elements of order $2$ and $\rho\notin \langle\sigma\rangle.$ 
Take $v=(1,-1,-1,1)\in V$ and set $S=\mathbb{F}_{7}v$. 
Then we know that $$\rho\cdot v=(-1,1,1,-1)=\sigma\cdot v=-v.$$ 
Hence, $S$ is a nonzero $\mathbb{F}_7\langle \rho, \sigma\rangle$-submodule of $V$. Therefore, both $$\mbox{Fix}_{S}(\rho)=0 \mbox{ and } \mbox{Fix}_{S}(\sigma)=0$$ since $2\neq 0$ in $\mathbb{F}_7$. By Proposition \ref{PPOO} (a3), $L(G)$ is inconsistent.}    
\end{example}

\begin{remark}\label{OORE}
\emph{In general, (a1) and (a2) in Proposition \ref{PPOO} do not imply (a3). For example, $L(A_4)$ in Example \ref{EX++} is consistent, and $V$ is simple and hence isotypic as an $\mathbb F_7A_4$-module. Thus both (a1) and (a2) are satisfied. However, (a3) fails.Therefore, in the nonabelian case, the consistency of $L(H)$ together with the isotypicity of $V$ does not, by itself, ensure the consistency of $L(G)$.  }
\end{remark}

In what follows, we first consider the case in which $H$ is abelian and establish a necessary and sufficient condition for the consistency of $L(G)$.
\begin{lemma}\label{OTH4.2}
Let $G=V\rtimes_{\theta} H$, where $H$ is a finite abelian group acting faithfully on $V$, $p\nmid |H|$ and $\emph{exp}(H)$ is square-free. Then $L(G)$ is consistent if and only if either $H=\{1\}$, or $H\cong C_q$ for some prime $q$ and $V$ is isotypic as an $\mathbb{F}_pH$-module.
\end{lemma}
\proof
Suppose first that $L(G)$ is consistent. By Proposition \ref{PRO4.1}, $V$ is isotypic as an $\mathbb{F}_pH$-module. Hence there exist a simple $\mathbb{F}_pH$-module $S\leq V$ and an integer $l\geq 1$ such that
$$V\cong S^{\oplus l}.$$ 
We first claim that $H$ acts
faithfully on $S$. Otherwise, there exists $1\neq h\in H$ that acts trivially on $S$. Since $V\cong S^{\oplus l}$ as an 
$\mathbb{F}_pH$-module, $h$ acts trivially on $V$, contrary to the fact that the action of $H$ on $V$ is faithful. 

Assume, for a contradiction, that $H\neq \{1\}$ and that $H$ is not cyclic of prime order. We claim that there exist two commuting elements $\rho,\sigma\in H$ of prime order such that $$\rho\notin\langle\sigma\rangle.$$ 
Indeed, if $|H|$ is divisible by 
two distinct primes $r$ and $s$, then choose elements $\rho, \sigma\in H$
of orders $r$ and $s$, respectively.  Since $H$ is abelian, $\rho$ and $\sigma$
commute, and clearly $\rho\notin\langle\sigma\rangle$.
On the other hand, if $|H|$ is divisible by only one prime $r$, then
$\mbox{exp}(H)=r$ since $\mbox{exp}(H)$ is square-free. As $H$ is a finite abelian group, it follows that $H\cong C^{t}_{r}$ 
for some $t\geq 1$. Since $H\not\cong C_r$, we have $t\geq 2$. Hence we may choose two elements $\rho,\sigma\in H$ of order $r$ such that
$\rho\notin\langle\sigma\rangle$. Note that $H$ is abelian. Therefore, $\rho$ and $\sigma$ commute.  Moreover, it is clear that $\mbox{Fix}_{S}(\rho)$ is an $\mathbb{F}_pH$-submodule of $S$ because $H$ is abelian. Note that $S$ is simple and the action of $H$ on $S$ is faithful. Then $\mbox{Fix}_{S}(\rho)\neq S$, and hence $\mbox{Fix}_{S}(\rho)=0$. Similarly, $\mbox{Fix}_{S}(\sigma)=0$. 
As $S$ is an $\mathbb{F}_pH$-submodule of $V$ and $\langle\rho,\sigma\rangle\leq H$, $S$ is also a nonzero 
$\mathbb{F}_p\langle \rho, \sigma\rangle$-submodule of $V$. Hence, (a3) in Proposition \ref{PPOO} is violated. Therefore, $H=\{1\}$ or $H\cong C_q$ for some prime $q$.

Conversely, if $H=\{1\}$, then $G\cong V$. Thus $L(G)$ is 
modular since $G$ is abelian, and therefore it is consistent. 
Now suppose that $H\cong C_q=\langle b\rangle$ for some prime $q$ 
and that $V$ is isotypic as an $\mathbb{F}_pH$-module. Put
$$\alpha=\theta(b)\in \mbox{GL}_{\mathbb{F}_p}(V).$$
As the action of $H$ on $V$ is faithful, $o(\alpha)=q$. Moreover, $p\nmid q$ and $q$ is square-free. Therefore, $L(G)$ is consistent by Theorem \ref{T001}.
\endproof

In view of Remark \ref{OORE}, it is natural to conjecture that if $G$ satisfies the assumptions of Proposition \ref{PPOO} and is consistent, then the finite group $H$ must be an abelian group.

\begin{proposition}\label{Propo+0}
Let $G=V\rtimes_{\theta}H$, where $H$ is a finite group acting faithfully on $V$ and
$p\nmid |H|$. If $H$ contains a subgroup $E\cong C_r\times C_r$ for some prime $r$, then $L(G)$ is inconsistent.
\end{proposition} 
\proof 
Since $E\leq H$, the restriction
$\theta|_E:E\longrightarrow \mbox{GL}_{\mathbb{F}_p}(V)$ 
defines a linear action of $E$ on $V$.
Thus $V$ is naturally an $\mathbb{F}_pE$-module.  
Moreover, since $p\nmid |H|$, we have $p\nmid |E|$. 
Hence, by Theorem \ref{cq1}, $V$ is semisimple as an $\mathbb{F}_pE$-module. 
Therefore, $$V=S_1\oplus\cdots\oplus S_t,$$ 
where each $S_i$ is a simple $\mathbb{F}_pE$-submodule of $V$. 
Note that the restricted action of $E$ on $V$ is faithful. 
Thus, there exists an integer $1\leq j\leq t$ such that 
the action of $E$ on $S_j$ is nontrivial.
Set $S=S_j$. Define 
$$\theta_S: E\longrightarrow \mbox{GL}_{\mathbb{F}_p}(S), \theta_S(e)=\theta(e)|_S$$
and let $$K=\mbox{ker}\theta_S=\{e\in E: e\cdot v=v \mbox{ for all } v\in S\}.$$

Let $F=\mathbb{F}_p[\theta_S(E)]$. We now prove that $F$ is a finite field. 
Since $E$ is abelian, the algebra $F$ is commutative. Set $0\neq f\in F$. As 
$f$ commutes with $\theta_S(e)$ for every $e\in E$, $\mbox{ker}f$ is an $\mathbb{F}_pE$-submodule of $S$.
We have $\mbox{ker}f=0$ since $f\neq 0$ and $S$ is simple. Hence, $f$ is invertible as a linear 
transformation on $S$ since $S$ is a finite-dimensional vector space. 
Note that $F$ is finite. There exists a positive integer $m$ such that $f^m=\mbox{id}_S$. Hence $f^{-1}=f^{m-1}\in F$.
Therefore, $F$ is a finite field. As the multiplicative group $F^{\times}$ of a finite field is cyclic and $\theta_S(E)\leq F^{\times}$, 
the image $\theta_S(E)$ is cyclic. 
Moreover, $E/K\cong \theta_S(E)$, $\theta_S(E)\neq \{1\}$ and 
$\mbox{exp}(\theta_S(E))\mid \mbox{exp}(E)=r$. 
Thus $\theta_S(E)\cong C_r$ and $|K|=r$, which means that 
$E\setminus K \neq \emptyset$. Hence, there exists $\sigma\in E\setminus K$.
Furthermore,  $E\setminus (K\cup \langle \sigma \rangle) \neq \emptyset$, and hence there also exists $\rho \in E\setminus (K\cup \langle \sigma \rangle)$.
Therefore, $\rho$ and $\sigma$ are commuting elements of order $r$ and
$\rho\notin\langle\sigma\rangle$. Thus $\langle\rho,\sigma\rangle=E$. On the other hand, as $E$ is abelian, $\mbox{Fix}_S(\rho)$ is an $\mathbb{F}_pE$-submodule of the simple module $S$. 
Hence $$\mbox{Fix}_{S}(\rho)= 0 \mbox{ or } \mbox{Fix}_{S}(\rho)=S.$$
If $\mbox{Fix}_{S}(\rho)=S$ then $\rho\in K$, a contradiction. 
Thus $\mbox{Fix}_{S}(\rho)= 0$. Similarly, $\mbox{Fix}_{S}(\sigma)=0$.
Therefore, by Proposition \ref{PRO4.2},  $L(G)$ is inconsistent. 
\endproof

\begin{theorem}[\cite{Rotman}]\label{Sylow}
Let $G$ be a finite group with $|G|=p^{i}m$, where $p$ is a prime and $p\nmid m$. Then: \\
(I) $G$ contains a subgroup of order $p^{i}$;\\
(II) Any two Sylow $p$-subgroups of $G$ are conjugate; \\
(III) $n_p\mid m$ and $n_p\equiv 1 (\emph{mod } p)$, where $n_p$ denotes the number of Sylow $p$-subgroups of $G$.\\
In particular, a Sylow $p$-subgroup is normal if and only if it is unique.
\end{theorem}

\begin{lemma}\label{LEMMA+1}
Let $G=V\rtimes_{\theta}H$, where $H$ is a finite group acting faithfully on $V$ and $p\nmid |H|$. If $H$ contains a nonabelian subgroup $B$ of order $rk$, where $r>k$ are distinct primes, then $L(G)$ is inconsistent.
\end{lemma}
\proof 
By Theorem \ref{Sylow}, $n_r\mid k$ and $n_r\equiv 1 (\mbox{mod } r)$ 
where $n_r$ denotes the number of Sylow $r$-subgroups of $B$. 
As $k<r$, it follows that $1\leq n_r\leq k<r$. Thus $n_r=1$.
Therefore, the unique Sylow $r$-subgroup $R$ of $B$ is normal.
Let $Q$ be a Sylow $k$-subgroup of $B$. Since $|R|=r, |Q|=k$ 
and $R\trianglelefteq B$, we have $|RQ|=|R||Q|=rk=|B|$. 
Therefore,
$$B=R\rtimes Q, R=\langle\rho\rangle\cong C_r \mbox{ and } Q=\langle\sigma\rangle\cong C_k,$$ 
where $\sigma\rho\sigma^{-1}=\rho^i$ for 
some integer $i$. Moreover, $i\not\equiv 1 (\mbox{mod } r)$ 
since $B$ is nonabelian. Set 
$$C=\mbox{Fix}_V(R)=\{v\in V: x\cdot v=v \mbox{ for all } x\in R\}.$$
Since $R\trianglelefteq B$, for $b\in B$, $x\in R$ and $v\in C$, 
$$x\cdot (b\cdot v)=b\cdot((b^{-1}xb)\cdot v)=b\cdot v.$$ 
Thus $b\cdot v\in C$, which means that $C$ is $B$-invariant. 
Hence $C$ is an $\mathbb{F}_pB$-submodule of $V$ since $C$ is a subspace of $V$. 
On the other hand, since $p\nmid |B|$, Theorem \ref{cq1} implies 
that $V$ is semisimple as an $\mathbb{F}_pB$-module. Note that $C$ 
is an $\mathbb{F}_pB$-submodule of $V$. Thus there exists an $\mathbb{F}_pB$-submodule $W$ such that 
$$V=C\oplus W.$$ Therefore, $\mbox{Fix}_W(R)=W\cap \mbox{Fix}_V(R)=W\cap C=0$.
Moreover, we claim that $W\neq 0$. Otherwise, $R$ 
would act trivially on $V$, contrary to the fact that 
the action of $H$ is faithful.

Now, we claim that $\sigma$ does not act trivially on $W$. 
Otherwise, $\sigma\cdot w=w$ for all $w\in W$. As $W$ is $B$-invariant, $\rho\cdot w\in W$. Hence, for every $w\in W$, 
$$\rho^i\cdot w
=(\sigma\rho\sigma^{-1})\cdot w
=\sigma\cdot(\rho\cdot(\sigma^{-1}\cdot w))
=\rho\cdot w$$ since $\sigma\rho\sigma^{-1}=\rho^i$. 
Thus $\rho^{i-1}\cdot w=w$ for all $w\in W$.
Note that $i\not\equiv 1 (\mbox{mod } r)$. 
Then $R=\langle\rho^{i-1}\rangle$, contrary to the 
fact that $\mbox{Fix}_W(R)=0$. 
Hence, there exists $y\in W$ such that $z=y-\sigma\cdot y\neq 0$. 

Put $$J=\langle(z,\sigma)\rangle, X=\{0\}\times R\mbox{ and }T=\langle X,J\rangle.$$
As $(z,\sigma)=(y,1)(0,\sigma)(-y,1)$, $(z,\sigma)$ is conjugate to $(0,\sigma)$. Thus $|J|=k$, which together with Proposition \ref{P02} implies that $J\in A(L(G))\subseteq J^+(L(G))$. Set $\tau=\sigma\rho\sigma^{-1}=\rho^i$. Then 
$o(\tau)=o(\rho)=r$, which means that $\tau\neq 1$ and $\langle\tau\rangle=R$. Moreover, we have 
$$[(z,\sigma),(0,\rho)]=(z-\tau\cdot z,\tau\rho^{-1})\in T.$$
Obviously, $(0,(\tau\rho^{-1})^{-1})\in T$ since $\tau\rho^{-1}\in R$.
Thus
\begin{equation}\label{LE00}
(z-\tau\cdot z,\tau\rho^{-1})(0,(\tau\rho^{-1})^{-1})=(z-\tau\cdot z,1)\in T\cap V. \end{equation}
Let $Y=\{w\in W: (w, 1)\in T\}=\{w\in W: (w, 1)\in T\cap V\}$. 
Since $V\trianglelefteq G$, we have
$T\cap V\trianglelefteq T$.
Moreover, $(0,\tau)\in X\leq T$. Hence, for every $(v,1)\in T\cap V$,
$$(0,\tau)(v,1)(0,\tau)^{-1}=(\tau\cdot v,1)\in T\cap V.$$
Therefore, the subspace $\{x\in V: (x,1)\in T\cap V\}$ is $\tau$-invariant. 
Note that $W$ is $B$-invariant and $\tau\in R\leq B$. Thus, $W$ is $\tau$-invariant.
Hence, $Y$ is also $\tau$-invariant. As $z\in W$, we have $z-\tau\cdot z\in W$. Therefore, by \eqref{LE00}, 
\begin{equation}\label{EQ-1}
z-\tau\cdot z \in Y.   
\end{equation}
Moreover, we know that $\mbox{Fix}_W(\tau)=\mbox{Fix}_W(R)=0$ since $\tau\neq 1$ and $\langle\tau\rangle=R$. 
We define the linear map $$f:W\longrightarrow W, f(w)=w-\tau\cdot w.$$
Thus $f$ is injective on $W$, which together with $Y\subseteq W$ and $Y$ is $\tau$-invariant implies that $f|_Y$ is injective on $Y$. Since $Y$ is 
finite-dimensional, $f|_Y$ is an automorphism of $Y$. From \eqref{EQ-1}, 
$f(z)=z-\tau \cdot z\in Y$. Since $f|_Y$ is surjective, 
there exists $u\in Y$ such that $$f|_Y(u)=f(u)=f(z).$$
As $u, z\in W$ and $f$ is injective on $W$, it follows that $u=z$. 
Therefore, $z\in Y$, and hence $(z, 1)\in T$. Thus, $(0,\sigma)=(-z,1)(z,\sigma)\in T$, which together with $X=\{0\}\times R\leq T$ implies that $$\{0\}\times B\leq T.$$

Let $U=\mbox{span}_{\mathbb{F}_p}\{b\cdot z: b\in B\}$. 
Since $\{0\}\times B\leq T$ and $(z,1)\in T\cap V$,  
$(0,b)(z,1)(0,b)^{-1}=(b\cdot z,1)\in T$ for every $b\in B$. 
Hence, $U\times \{1\}\leq T$, which implies that $U\rtimes B\leq T$. 
The reverse inclusion follows from
$T=\langle(0,\rho),(z,\sigma)\rangle$.
Therefore, $T=U\rtimes B$. We denote $M_1=\{0\}\times B$ and $M_2=U\rtimes R$.
Then $X\leq M_1, M_2<T$ because $(z, 1)\in T\setminus M_1$ and $(0,\sigma)\in T\setminus M_2$. As $T=U\rtimes B$, $M_1\vee M_2=T$. Thus $T=X\vee J$ is 
not join-irreducible in the interval $[X,G]$. However, $J\in J^+(L(G))$. Therefore, $L(G)$ is inconsistent.
\endproof

By Lemma \ref{LEMMA+1} and Proposition \ref{Propo+0}, we have the following lemma.
\begin{lemma}\label{LP}
Let $G=V\rtimes_{\theta}H$, where $H$ is a finite group acting faithfully on $V$, $p\nmid |H|$ and $\emph{exp}(H)$ is square-free. If $L(G)$ is consistent then $H$ is abelian.
\end{lemma}
\proof
Assume that $H$ is nonabelian. If $H$ contains a subgroup $E\cong C_r\times C_r$ for some prime $r$, then Proposition \ref{Propo+0} implies that $L(G)$ is inconsistent, 
a contradiction. Now, suppose that $H$ contains no elementary abelian subgroup of order $r^2$ for any prime $r$. Since $H$ is nonabelian 
and has square-free exponent, Case (a) of Lemma 3 in Jabara \cite{Jabara} implies that $H$ contains a nonabelian subgroup $B$ of order $rk$, where $r>k$ are distinct primes. 
By Lemma \ref{LEMMA+1}, $L(G)$ is inconsistent, a contradiction. Therefore, $H$ is abelian.
\endproof

\begin{lemma}\label{TH4.2}
Let $G=V\rtimes_{\theta} H$, where $H$ is a finite group acting faithfully on $V$, $p\nmid |H|$ and $\emph{exp}(H)$ is square-free. Then $L(G)$ is consistent if and only if either $H=\{1\}$, or $H\cong C_q$ for some prime $q$ and $V$ is isotypic as an $\mathbb{F}_pH$-module.
\end{lemma}
\proof Suppose first that $L(G)$ is consistent. 
By Lemma \ref{LP}, $H$ is abelian. Lemma \ref{OTH4.2} then implies that 
either $H=\{1\}$, or $H\cong C_q$ for some prime $q$ and $V$ is isotypic as an 
$\mathbb{F}_pH$-module.

 Conversely, if either of these conditions holds, then $H$ 
is abelian. Thus, Lemma \ref{OTH4.2} shows that $L(G)$ is consistent.
\endproof

\begin{remark}
\emph{Let $G$ be a finite group and $L(G)$ be atomistic. Then 
every nontrivial subgroup $K\leq G$ admits an irredundant 
generating expression
$$K=\langle P_1,\cdots,P_r\rangle, P_i\in A(L(G)),$$
where $$A(L(G))=\{P\leq G: |P| \mbox{ is prime}\}.$$
Reuter \cite{Reuter} states that a finite lattice $L$ is 
consistent if and only if it satisfies the Kurosh-Ore 
replacement property. Therefore, $L(G)$ is consistent if and only if 
for any two irredundant generating expressions of $K\leq G$,
$$K=\langle P_1,\cdots,P_r\rangle=\langle Q_1,\cdots,Q_s\rangle, P_i, Q_j\in A(L(G)),$$ 
and for every $P_i$ there exists $Q_j$ such that
$$K=\langle P_1,\cdots,P_{i-1},Q_j, P_{i+1},\cdots,P_r
\rangle .$$
Consequently, Theorem \ref{T001} and Lemma \ref{TH4.2} provide structural characterizations of groups satisfying the Kurosh-Ore replacement property.} 
\end{remark}

\begin{lemma}\label{RE++1}
Let $G=V\rtimes_{\theta}H$, where $H$ is a finite group acting faithfully on $V$ and $p\nmid |H|$. If $H$ contains 
a cyclic subgroup of order $q^2$ for some prime $q$, then $L(G)$ is inconsistent.
\end{lemma}
\proof
Suppose that $H$ contains a cyclic subgroup $K=\langle c \rangle$ of order $q^2$. Put $a=c^q$. Then $o(a)=q$. 
As $p\nmid |H|$, we have $p\nmid |K|$. By Theorem \ref{cq1}, $V$ is semisimple as an $\mathbb{F}_pK$-module. 
Note that $H$ acts faithfully on $V$ and $a\neq 1$. Thus, there exists a simple $\mathbb{F}_pK$-submodule $S\leq V$ 
such that $a$ acts nontrivially on $S$. We first claim that $$\mbox{Fix}_S(a)=0.$$ Indeed, since $K$ is abelian, $\mbox{Fix}_S(a)$ is an 
$\mathbb{F}_pK$-submodule of $S$. Since $S$ is simple, $\mbox{Fix}_S(a)=0$ or $\mbox{Fix}_S(a)=S$. 
However, $a$ acts nontrivially on $S$. Hence, $\mbox{Fix}_S(a)\neq S$, which means that $\mbox{Fix}_S(a)=0$.

Let $0\neq y\in S$ and put $x=y-c\cdot y$. Then we have $x\neq 0$. Otherwise, $c\cdot y=y$, 
which would imply $a\cdot y=c^q \cdot y=y$, a contradiction. Denote $$J=\langle (x, c)\rangle, U=\langle (0, a)\rangle \mbox{ and }T=\langle U, J\rangle.$$
As $(x, c)=(y, 1)(0, c)(-y,1)$, the element $(x,c)$ is conjugate to $(0,c)$. Thus, $J$ is cyclic of order $q^2$. By Proposition \ref{P02}, $J\in J^{+}(L(G))$. 
We now prove that $$T=S\rtimes_{\theta} K.$$
As $S$ is an $\mathbb{F}_pK$-submodule of $V$, we have $x=y-c\cdot y\in S$. Moreover, $a, c\in K$, so both $(0,a)$ and $(x,c)$ belong to $S\rtimes_{\theta} K$. 
Hence, $T\leq S\rtimes_{\theta} K$. On the other hand, we have 
$$\begin{aligned}
[(0,a),(x,c)]
&=(0,a)(x,c)(0,a)^{-1}(x,c)^{-1}\\
&=(a\cdot x,ac)(0,a)^{-1}(x,c)^{-1}\\
&=(a\cdot x,c)(-c^{-1}\cdot x,c^{-1}) \mbox{ since } ac=ca\\
&=(a\cdot x-x,1).
\end{aligned}$$

Define $R=\{z: (z,1)\in T\cap (S\times\{1\})\}$. Since $x\neq 0$ and $\mbox{Fix}_S(a)=0$, we have $a\cdot x-x\neq 0$. Therefore, $R\neq 0$. 
Now consider the natural projection $$\pi: S\rtimes_{\theta} K \longrightarrow K, \pi(s, k)=k.$$
Since $(x, c)\in T$ and $K=\langle c \rangle$, we have $\pi(T)=K$. We next show that $R$ is an $\mathbb{F}_pK$-submodule of $S$.
Note that $S$ is an elementary abelian $p$-group. The subgroup $R$ is an $\mathbb{F}_p$-subspace of $S$. 
Moreover, as $S\times \{1\} \unlhd S\rtimes_\theta K$, $T\cap (S\times \{1\})\unlhd T$. For any $k\in K$, $\pi(T)=K$ implies that there 
exists some $v\in S$ such that $(v, k)\in T$. Thus, for every $r\in R$, $$(v, k)(r, 1)(v, k)^{-1}=(k\cdot r, 1)\in T\cap (S\times\{1\}),$$
which means that $k\cdot r\in R$. Therefore, $R$ is $K$-invariant, and hence 
\begin{equation}\label{RRR}
R \mbox{ is }\mathbb{F}_pK\mbox{-submodule of } S.
\end{equation} Since $S$ is simple and $R\neq 0$, we have $R=S$. 
Hence $S\times\{1\}\leq T$, which together with $\pi(T)=K$ implies that $S\rtimes_{\theta} K\leq T$. Note that $T\leq S\rtimes_{\theta} K$. 
Consequently, $T=S\rtimes_{\theta} K$.

Let $A=\{0\}\times K$ and $B=S\rtimes_\theta \langle a\rangle$. We show that $$A\prec T \mbox{ and }B\prec T.$$ 
As $$\frac{|T|}{|B|}=\frac{|S|\cdot |K|}{|S|\cdot|\langle a\rangle|}=\frac{|K|}{|\langle a\rangle|}=q,$$ $B$ is maximal in $T$, which means that $B\prec T$. 
Suppose that $A\leq M\leq T$ for some subgroup $M$. Denote $W=\{x: (x,1)\in M\cap (S\times\{1\})\}$. Let $(s,k)\in M$. Since $A=\{0\}\times K\le M$, we have $(0,k)\in M$. Hence
$(s,1)=(s,k)(0,k)^{-1}\in M\cap (S\times\{1\})$, which means that $s\in W$. Therefore, $M=W\rtimes_\theta K$. Furthermore, similar to the proof of \eqref{RRR}, $W$ is an $\mathbb{F}_pK$-submodule of $S$. As $S$ is 
simple, $W=0$ or $W=S$. Hence, $M=A$ or $M=T$. Therefore, $A\prec T$. 

Obviously, $A\neq B$ and $A\cap B=A\wedge B=\{0\}\times \langle a\rangle=U$. Moreover, $T\geq A\vee B\geq S\rtimes_{\theta} K=T$ implies $T=A\vee B$. 
Hence, $A\wedge B=U, A\vee B=T, A, B\prec T$ and $U\vee J=T$. Consequently, by Theorem \ref{TT01}, $L(G)$ is inconsistent.
\endproof

\quad\\
\textbf{Proof of Theorem \ref{TH4.3}:}
Suppose first that $L(G)$ is consistent. We claim that $\mbox{exp}(H)$ is square-free. Otherwise, there exists a prime $q$ such that $q^2\mid\mbox{exp}(H)$. 
Hence $H$ contains a cyclic subgroup of order $q^2$. By Lemma \ref{RE++1}, $L(G)$ is inconsistent, a contradiction. Therefore, $\mbox{exp}(H)$ is square-free.
Using Lemma \ref{TH4.2}, we conclude that either $H=\{1\}$, or $H\cong C_q$ for some prime $q$ and $V$ is isotypic as an $\mathbb{F}_pH$-module.

Conversely, suppose that either $H=\{1\}$, or $H\cong C_q$ for some prime $q$ and $V$ is isotypic as an $\mathbb{F}_pH$-module. In either case, $\mbox{exp}(H)$ is square-free. 
Hence, by the characterization of consistency in Lemma \ref{TH4.2}, $L(G)$ is consistent.
\endproof

We conclude this subsection by investigating the relation between consistency 
of subgroup lattices and Frobenius groups. We first recall a known criterion for Frobenius groups. 
For more related material, see \cite{KKK, Scott}

\begin{theorem}[\cite{KKK}]\label{Frob}
Let $G$ be the finite semidirect product of a nontrivial subgroup $H$ with
a nontrivial normal subgroup $K$. Then the following statements are equivalent:\\
(i) $G$ is a Frobenius group with Frobenius complement $H$ and Frobenius
kernel $K$.\\
(ii) $\emph{Fix}_K(h)=\{1_K\}$ for all $1\neq h\in H$.
\end{theorem}

Then Theorems \ref{TH4.3} and \ref{Frob} yield the following corollary.
\begin{corollary}\label{corr4.2}
Let $G=V\rtimes_\theta H$ satisfy the hypotheses of Theorem \ref{TH4.3}, and suppose that $H\neq\{1\}$.
If $L(G)$ is consistent, then $G$ is a Frobenius group with Frobenius kernel $V$ and Frobenius complement $H$.
\end{corollary}
\proof
By Theorem \ref{TH4.3}, $H\cong C_q$ for some prime $q$, and $V$ is isotypic as an $\mathbb{F}_pH$-module. Thus $V\cong S^{\oplus t}$ for 
some simple $\mathbb{F}_pH$-module $S$ and some positive integer $t$. We claim that $$\mbox{Fix}_S(h)=0$$ for every $1\neq h\in H$. 
Indeed, by Theorem \ref{TH4.3}, $H=\langle h\rangle$. It follows that $\mbox{Fix}_S(h)$ is an $\mathbb{F}_pH$-submodule of $S$. Then $\mbox{Fix}_S(h)=0$ or $\mbox{Fix}_S(h)=S$ since $S$ is simple. If $\mbox{Fix}_S(h)=S$, then the action of $H$ on $S$ is not faithful, 
and hence the action of $H$ on $V$ is also not faithful, a contradiction. Hence $\mbox{Fix}_S(h)=0$, which implies that 
$\mbox{Fix}_V(h)=0$ for every $1\neq h\in H$. 
Consequently, by Theorem \ref{Frob}, $G$ is a Frobenius group with Frobenius complement $H$ and Frobenius kernel $V$.
\endproof

\begin{remark}\label{REmark}
\emph{The converse of Corollary \ref{corr4.2} does not hold in general. For example, let $V=\mathbb{F}^2_{7}$ and $H=C_3=\langle b\rangle$. 
Define $\alpha\in \mbox{GL}_{\mathbb{F}_{7}}(V)$ by
$$\alpha=\begin{bmatrix}
	    2 & 0 \\
	    0 & 4
	\end{bmatrix}.$$ 
Then $o(\alpha)=3$ and $\alpha^2=\begin{bmatrix}
	    4 & 0 \\
	    0 & 2
	\end{bmatrix}.$
Hence, the action of $H$ on $V$ is faithful. Furthermore, it is easy to verify that $\mbox{Fix}_V(b)=\mbox{Fix}_V(b^2)=0$. 
By Theorem \ref{Frob}, $G=V\rtimes_{\alpha} H$ is a Frobenius group. However, the minimal polynomial $m_{\alpha}(x)=(x-4)(x-2)$ is 
reducible over $\mathbb{F}_7$.
By Theorem \ref{T001}, $L(G)$ is inconsistent. Therefore, the converse of Corollary \ref{corr4.2} does not hold in general.
The above example also shows that the Frobenius property does not imply the isotypicity of $V$ as an $\mathbb{F}_pH$-module by Theorem \ref{TH4.3}.}
\end{remark}

\section{Conclusions}
In this article, we established several equivalent characterizations of consistency 
for lattices of finite length via forbidden sublattices 
and CEP associated with the canonical 
closure operator. We then applied these results to semidirect products 
$G=V\rtimes H$ and obtained a complete characterization of the consistency 
of subgroup lattices in terms of the structure of the acting group $H$ 
and the $\mathbb{F}_pH$-module $V$. These results deepen the connection 
between the consistency of subgroup lattices and the representation 
theory of finite groups.
%\section*{Acknowledgments}
%This work was partially supported by  the National Natural Science Foundation of China (Nos. 12471440 and 12671562).

\section*{AI Declaration}
The authors declare that no generative AI or AI-assisted technologies were used in the preparation of this manuscript.
%\section*{ORCID}
%\quad\\
%Peng He:  https://orcid.org/0009-0003-6101-5644\\
%Xue-ping Wang: https://orcid.org/0000-0003-4376-9859

\end{document}